\documentclass{amsart}

\usepackage{graphicx}
\usepackage{tikz-cd}
\usepackage[latin1]{inputenc}

\newtheorem{thm}{Theorem}[section]
\newtheorem{cor}[thm]{Corollary}
\newtheorem{lem}[thm]{Lemma}
\newtheorem{prop}[thm]{Proposition}
\newtheorem{exa}[thm]{Example}
\theoremstyle{definition}
\newtheorem{dfn}[thm]{Definition}
\theoremstyle{remark}
\newtheorem{rem}[thm]{Remark}
\numberwithin{equation}{section}

\begin{document}

\title[]{Liftings of knots in $S_{g} \times S^{1}$ and covering of virtual knots}%
\author{S.Kim}

\address{S.Kim, Jilin university}%
\email{kimseongjeong@jlu.edu.cn}%


\maketitle

\begin{abstract}
A virtual link diagram is called {\em (mod $m$) almost classical} if it admits a (mod $m$) Alexander numbering. In \cite{BodenGaudreauHarperNicasWhite}, it is shown that Alexander polynomial for almost classical links can be defined by using the homology of the associated infinite cyclic cover. On the other hand, in \cite{NaokoKamada} an infinite family of $m$ fold covering over a virtual knot is constructed so that it is mod $m$ almost classical link for all $m$ by using oriented cut point. In this paper, we study liftings of knots in $S_{g} \times S^{1}$. Liftings of knots in $S_{g} \times S^{1}$ provide another way to obtain a family of $m$-fold coverings over a given virtual knots.\\[2mm]

\end{abstract}

\section{Introduction}

A virtual link diagram is called {\em (mod $m$) almost classical} if it admits a (mod $m$) Alexander numbering. In \cite{BodenGaudreauHarperNicasWhite}, it is shown that Alexander polynomial for almost classical links can be defined by using the homology of
the associated infinite cyclic cover. In particular, it satisfies a skein relation. Similarly one can construct an Alexander polynomial for mod $m$ almost classical links.

On the other hand, in \cite{Dye1,Dye2} the notions of an oriented cut point and a cut system for a virtual link diagram, which is an extension of (unorieted) cut points, are introduced. In \cite{NaokoKamada} it is proved that any two cut systems for a given virtual link diagram are equivalent up to finitely many local moves and an infinite family of $m$-fold covering over a virtual knot is constructed so that it is mod $m$ almost classical link for all $m$.

 We are interested in links in $S_{g}\times S^{1}$, where $S_{g}$ is an oriented surface of genus $g$. In \cite{Kim}, by the author diagrams on a plane and local moves for links in $S_{g}\times S^{1}$, which present equivalence relations for links in $S_{g}\times S^{1}$, are constructed.

This paper is contributed to construct another way to obtain a family of $m$-fold coverings over a given virtual knots, which are mod $m$ almost classical, by using knots in $S_{g} \times S^{1}$.

The present paper is organized as follows: In Section 2, we introduce the notion of almost classical links and oriented cut points. In Section 3, the basic notions and properties of links in $S_{g}\times S^{1}$ are introduced. In Section 4, the covering over knots in $S_{g} \times S^{1}$ is constructed by using heights of arcs of them. In Section 5, a map from virtual knots to knots in $S_{g}\times S^{1}$ by using oriented cut points. We will show that the $m$-fold covering over knots in $S_{g} \times S^{1}$ of degree $0$ constructed in Section 4 is almost classical.

\section{Almost classical links and cut system}

In this paper, arcs between classical crossings are called {\em arcs}. This part follows the notions in \cite{BodenGaudreauHarperNicasWhite}

 \subsection{Almost classical links}
\begin{dfn}
An {\em Alexander numbering} (respectively, {\em mod $m$ Alexander numbering}) of an oriented virtual knot diagram $D$ is an assignment of each arc of $D$ to a number in $\mathbb{Z}$ (respectively, in $\mathbb{Z}_{m}$) so that the numbers of four arcs around each classical crossing are assigned as shown in Fig.~\ref{fig:rule-Alex-num} for some $i\in \mathbb{Z}$ (respectively, $i \in \mathbb{Z}_{m}$). 

  \begin{figure}[h!]
\begin{center}
 \includegraphics[width = 3.5cm]{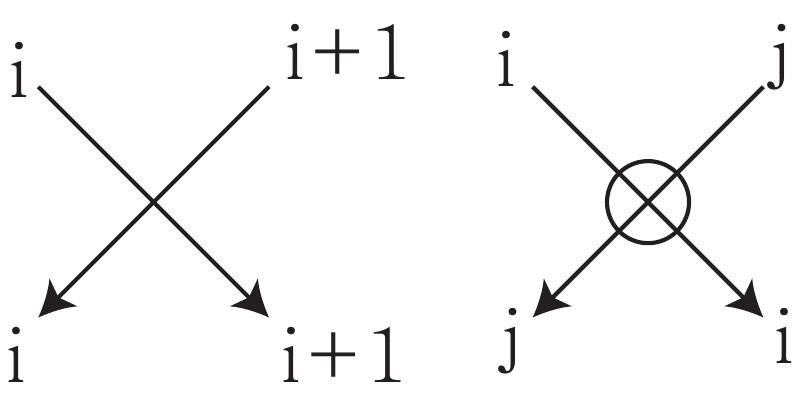}
\end{center}
\caption{Rules for the Alexander numbering on a virtual knot diagram}\label{fig:rule-Alex-num}
\end{figure}

If a diagram admits Alexander numbering (respectively, {\em mod $m$ Alexander numbering}), then it is called {\em almost classical diagram} (respectively, {\em mod $m$ almost classical diagram})
\end{dfn}

\begin{exa}
It is known that the left diagram in Fig.~\ref{fig:exa-Alexander} does not admit Alexander numbering.
The right diagram in Fig.~\ref{fig:exa-Alexander} admits an Alexander numbering as illustrated and hence it is almost classical.

  \begin{figure}[h!]
\begin{center}
 \includegraphics[width = 5cm]{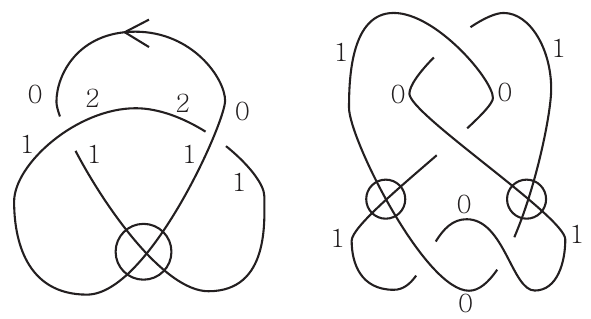}
\end{center}
\caption{Non-Alexander numberable and almost classical diagrams}\label{fig:exa-Alexander}
\end{figure}
\end{exa}

\begin{exa}
It is known that the virtual knot diagram in Fig.~\ref{fig:exa-mod3-Alex-num} is not almost classical diagram. But one can find a mod $3$ Alexander numbering as illustrated in Fig.~\ref{fig:exa-mod3-Alex-num} and hence it is mod $3$ almost classical diagram. 

  \begin{figure}[h!]
\begin{center}
 \includegraphics[width = 5cm]{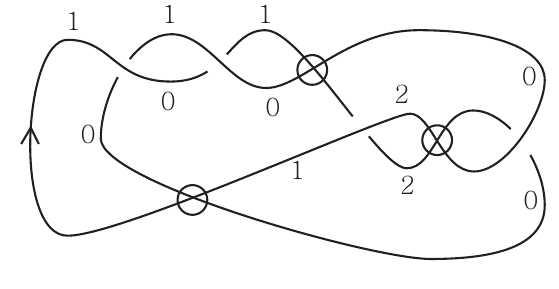}
\end{center}
\caption{Mod $3$ almost classical diagram}\label{fig:exa-mod3-Alex-num}
\end{figure}
\end{exa}

\begin{dfn}
A virtual link $L$ is {\em almost classical} (respectively, {\em mod $m$ almost classical}) if there is an almost classical (respectively, mod $m$ almost classical) virtual link diagram of $L$.
\end{dfn}

It is clear that any almost classical virtual link diagram is mod $m$ almost classical. One can show that a virtual link diagram is checkerboard colorable if and only if it is mod $2$ almost classical. 

\begin{prop}[\cite{BodenGaudreauHarperNicasWhite}]
For a virtual knot or link $K$, the followings are equivalent
\begin{enumerate}
\item $K$ is almost classical
\item $K$ is homological trivial as a knot or a link in $\Sigma \times [0,1]$.
\item $K$ is the boundary of a connected, oriented surface $F$ embedded in $\Sigma \times [0,1]$, where $\Sigma$ is the underlying surface associated to an Alexander numberable diagram of $K$.

\end{enumerate}
The surface $F$ is called {\em a Seifert surface for $K$}.
\end{prop}

\begin{prop}[\cite{BodenGaudreauHarperNicasWhite}]
A knot $K$ in a thickened surface $\Sigma \times [0,1]$ is mod $p$ almost classical if and only if it is homologically trivial as an element in $H_{1}(\Sigma; \mathbb{Z}_{p})$.
\end{prop}

In \cite{BodenGaudreauHarperNicasWhite} the Alexander polynomial $\Delta_{K}(t)$ is defined by using Seifert surfaces for almost classical links embedded in the thickened surface $\Sigma \times [0,1]$. In particular, it satisfies the usual skein relation $\Delta_{K_{+}}(t) -\Delta_{K_{-}}(t) = (t-t^{-1})\Delta_{K_{0}}(t)$.

\subsection{Cut system and $m$-fold cyclic covering}
In general, virtual link diagrams are not almost classical. In \cite{Dye1, Dye2}, the notion of cut points is introduced by H. Dye and it is possible to do an Alexander numbering on virtual link diagrams with cut points.
\begin{dfn}
An {\em oriented cut point} or simply {\em cut point} is an arrow on an arc of a virtual knot diagram which indicates a local orientation of the arc as illustrated in Fig.~\ref{fig:cutpoint}. An oriented cut point is called {\em coherent} (respectively, {\em incoherent}) if the local orientation indicated by the cut point agrees on (respectively, disagree on) the orientation of the virtual link diagram.
 \begin{figure}[h!]
\begin{center}
 \includegraphics[width = 4cm]{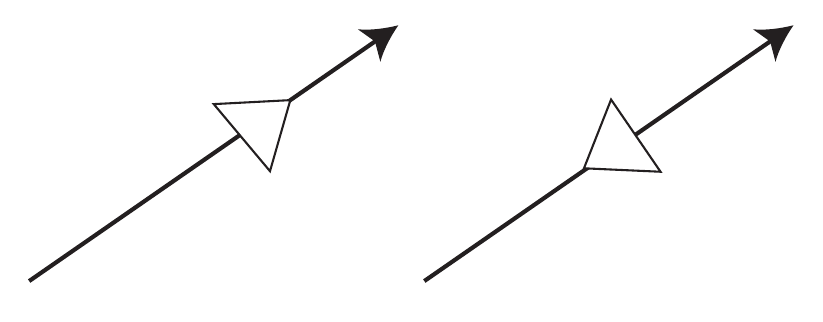}
\end{center}
\caption{Oriented cut points}\label{fig:cutpoint}
\end{figure}
\end{dfn}

\begin{dfn}
Let $D$ be a virtual link diagram and $P$ a set of oriented cut points of $D$. The set $P$ of cut points is called a {\em cut system} if $D$ admits an Alexander numbering such that at each oriented cut point, the number increases by one in the direction of the oriented cut point as described in Fig. \ref{fig:rule-Alex-num-cutpt}. Such an Alexander numbering is called {\em an Alexander numbering of a virtual link diagram with a cut system}. 
 \begin{figure}[h!]
\begin{center}
 \includegraphics[width = 6cm]{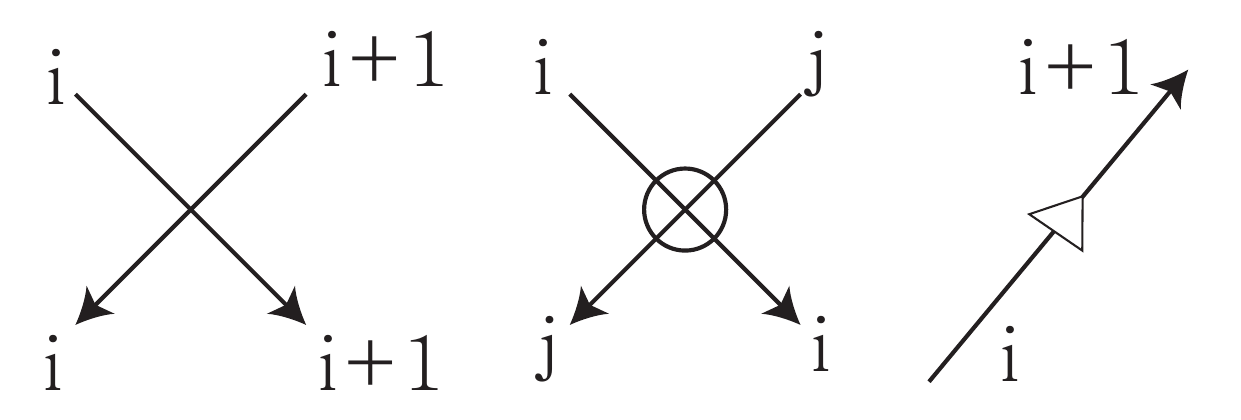}
\end{center}
\caption{Oriented cut points}\label{fig:rule-Alex-num-cutpt}
\end{figure}
\end{dfn}

\begin{dfn}
{\em A standard cut system} of a virtual link diagram is a cut system which is obtained by putting coherent and incoherent oriented cut points as described in the left of Fig.~\ref{fig:standard-cut-system} around each virtual crossing. It provides a cut system of the diagram and an Alexander numbering is given as described in the right of Fig.~\ref{fig:standard-cut-system} around each virtual crossing.
\begin{figure}[h!]
\begin{center}
 \includegraphics[width = 5cm]{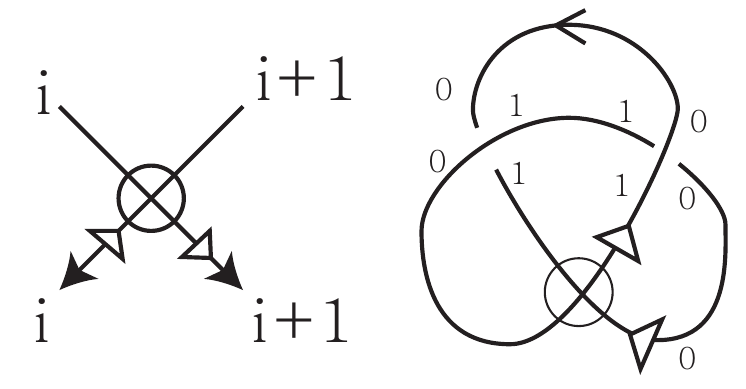}
\end{center}
\caption{Standard cut system}\label{fig:standard-cut-system}
\end{figure}
\end{dfn}

\begin{rem}
The idea of standard cut system is to assign an Alexander numbering around virtual crossings considering it as classical crossings. Note that every classical link diagram is almost classical, and hence, a virtual knot diagram with the standard cut system admits an Alexander numbering.
\end{rem}

\begin{thm}[\cite{NaokoKamada}]
Two cut systems of the same virtual link diagram are related by a sequence of oriented cut point moves described in Fig. \ref{fig:cutpt-moves}.

\begin{figure}[h!]
\begin{center}
 \includegraphics[width = 10cm]{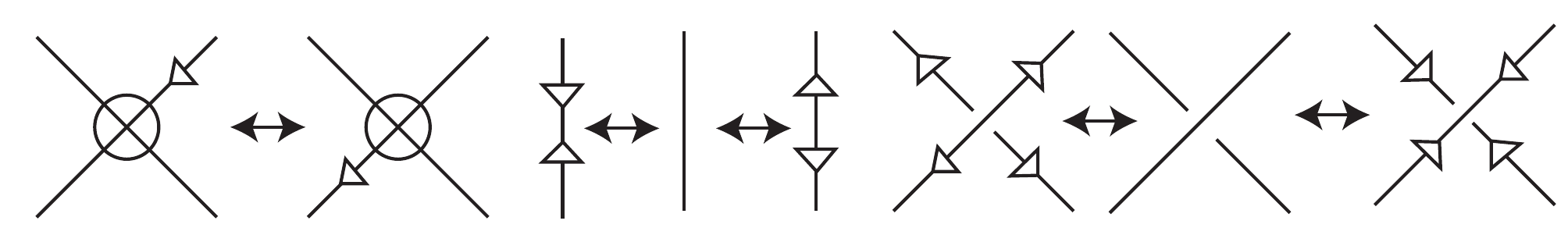}
\end{center}
\caption{Oriented cut point moves}\label{fig:cutpt-moves}
\end{figure}
\end{thm}

From the above theorem it follows that every cut system of a virtual link diagram is related with the standard cut system by a sequence of oriented cut point moves. 
\begin{cor}
Let $D$ be a virtual link diagram and $P$ a cut system of $D$. The number of coherent cut points of $P$ equals that of incoherent cut points of $P$.
\end{cor}

By N. Kamada an m-fold cyclic covering (virtual link) diagram $\phi_{m}(D,P)$ of an oriented virtual knot diagram with a cut system $P$ is constructed. The following proposition is proved, for details, see \cite{NaokoKamada}.

\begin{prop}
For a virtual link diagram $D$ with a cut system $P$, an $m$-fold cyclic covering virtual link diagram $\phi_{m}(D,P)$ is mod $m$ almost classical.
\end{prop}

\section{Basic notions for links in $S_{g} \times S^{1}$}

\subsection{Links in $S_{g} \times S^{1}$ and their diagrams}

\begin{dfn}
Let $S_{g}$ be an oriented surface of genus $g$. {\em A link $L$  in $S_{g} \times S^{1}$} is a pair of $S_{g} \times S^{1}$ and a smooth embedding of a disjoint union of $S^{1}$ into $S_{g} \times S^{1}$. We denote it by $(L, S_{g} \times S^{1})$. Each image of $S^{1}$ in $S_{g} \times S^{1}$ is called {\em a component} of $L$. A link of one component is called {\em a knot in $S_{g} \times S^{1}$.}
\end{dfn}

\begin{dfn}
Let $(L, S_{g} \times S^{1})$ and $(L', S_{g'} \times S^{1})$ be two links in $S_{g} \times S^{1}$ and $S_{g'}\times S^{1}$, respectively. We call $(L, S_{g} \times S^{1})$ and $(L', S_{g'} \times S^{1})$ {\em equivalent,} if $(L', S_{g'} \times S^{1})$ can be obtained from $(L, S_{g} \times S^{1})$ by a series of ambient isotopies and stabilization/destabilization of $S_{g} \times S^{1}$.
\end{dfn}


By the {\em destabilization for $(L, S_{g} \times S^{1})$} we mean the following:
Let $C$ be a non-contractible simple closed circle on the surface $S_{g}$ such that there exists a torus $T$ in $S_{g} \times S^{1}$ homotopic to the torus $C \times S^{1}$ and not intersecting the link $L$. Then our destabilization is cutting of the manifold $S_{g} \times S^{1}$ along the torus $C \times S^{1}$ and pasting of two newborn components of boundary by $D^{2} \times S^{1}$.
The {\em stabilization for $S_{g} \times S^{1}$} is converse operation to the destabilization.

First let us construct diagrams for $(L, S_{g}\times S^{1})$ on the surface $S_{g}$ as follows: Let $L$ be an (oriented) link in $S_{g} \times S^{1}$. Assume that an orientation is given on $S^{1}$. Suppose that $x_{0} \in S^{1}$ is a point such that $S_{g} \times \{x_{0}\} \cap L$ is a set of finite points with no transversal points to $S_{g} \times \{x_{0}\}$. 
\begin{figure}[h]
\begin{center}
 \includegraphics[width = 12cm]{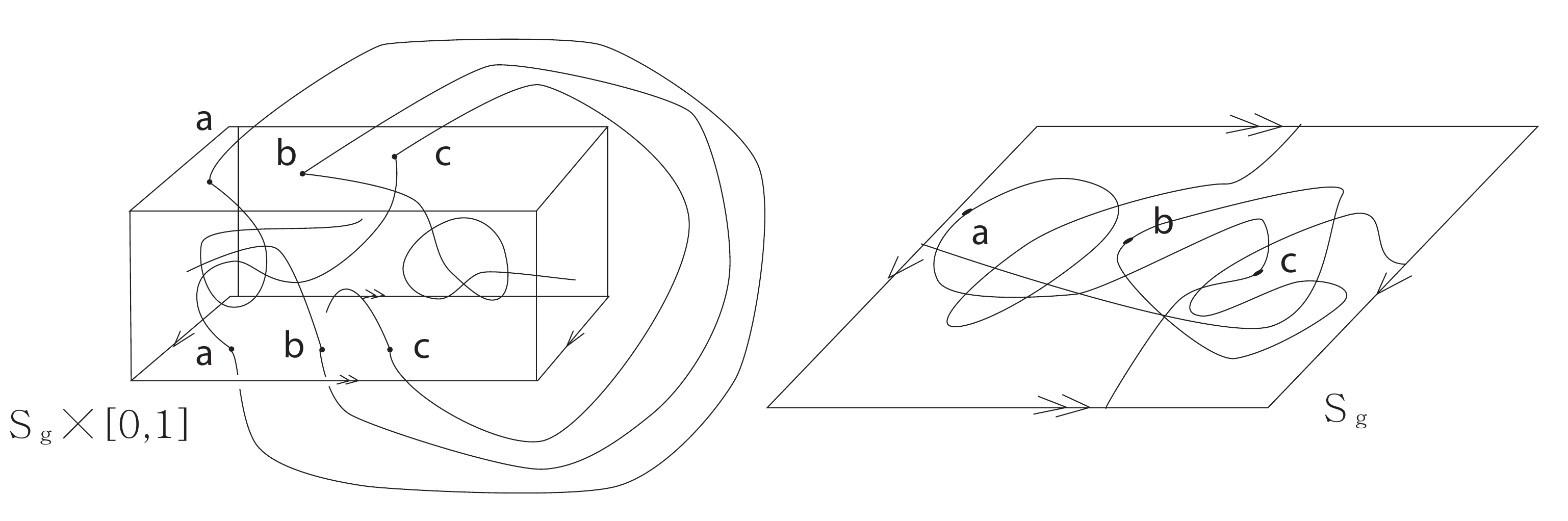}

\end{center}

\caption{Schematic figures of links in $S_{g} \times S^{1}$ and its projection on the surface $S_{g}\times \{0\}$}\label{fig:SS-diag-rel-1}
\end{figure}
Then there exists a natural diffeomorphism $f$\label{def:cutting-map} from $(S_{g} \times S^{1} )\backslash$  $(S_{g} \times \{x_{0}\})$ to $S_{g} \times (0,1) \subset S_{g} \times [0,1]$. Let $M_{L}$ $=$ $\overline{f((S_{g} \times S^{1}) - (S_{g} \times \{x_{0}\}))} \cong S_{g} \times [0,1]$. Then $\overline{f(L)}$ in $M_{L}$ consists of finitely many circles and arcs with two boundaries on $S_{g} \times \{0\}$ and $S_{g} \times \{1\}$.

Let $D_{\overline{f(L)}}$ be the image of a projection of $\overline{f(L)}$ onto the $S_{g} \times \{0\}$. Note that two boundary points $(p,0)$ and $(p,1)$ for some $p\in S_{g}$ of arcs of $\overline{f(L)}$ are projected to the same point, we call it {\em a vertex}. It follows that the diagram $D_{\overline{f(L)}}$ of $L$ on $S_{g}$ has $n$ circles with vertices corresponding to two boundary points on $S_{g} \times \{0\}$ and $S_{g} \times \{1\}$ and $m$ circles as described in the right of Fig.~\ref{fig:SS-diag-rel-1}.

\begin{figure}[h!]
\begin{center}
 \includegraphics[width = 8cm]{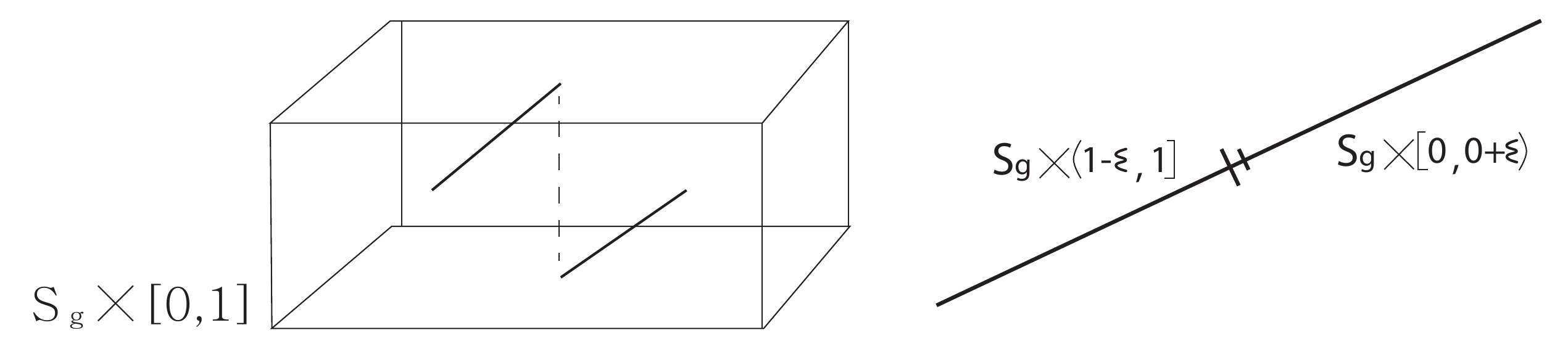}

\end{center}
\caption{Local image near to an intersection of a link with $S_{g} \times \{x_{0}\}$ and a corresponding double line}\label{fig:Vertex}
\end{figure}

Note that two arcs near to a vertex are the images of arcs near to $S_{g} \times \{0\}$ and $S_{g} \times \{1\}$ in $M_{L} \cong S_{g} \times [0,1]$, respectively, as described in Fig.~\ref{fig:Vertex}. We change the vertex to two short lines, which we call  {\it a double line}, so that if one of the lines is connected with an arc which is near to $S_{g} \times \{1\}$, then the line is longer than the another, as describe in Fig.~\ref{fig:Vertex}.

Since $D_{\overline{f(L)}}$ is a framed 4-valent graph with double lines on the surface $S_{g}$, which comes from $\overline{f(L)}$ in $S_{g} \times [0,1]$, one can give over/under information to each 4-valent vertex. That is, a link $L$ in $M_{L}$ has a knot diagram with double lines on $S_{g}$. Simply we call it {\em a diagram on $S_{g}$ with double lines.} 

\begin{prop}[M. K. Dabkowski, M. Mroczkowski (2009) \cite{DabkowskiMroczkowski}, Kim (2018) \cite{Kim}]\label{thm:diag_on_surface}
   Let $(L, S_{g} \times S^{1})$ and $(L', S_{g} \times S^{1})$ be two links in $S_{g} \times S^{1}$. Let $D_{L}$ and $D_{L'}$ be diagrams of $L$ and $L'$ on $S_{g}$, respectively. Then $(L, S_{g} \times S^{1})$ and $(L', S_{g} \times S^{1})$ are isotopic if and only if $D_{L'}$ can be obtained from $D_{L}$ by applying finitely many moves in Fig.~\ref{moves1}.

  \begin{figure}[h!]
\begin{center}
 \includegraphics[width = 9cm]{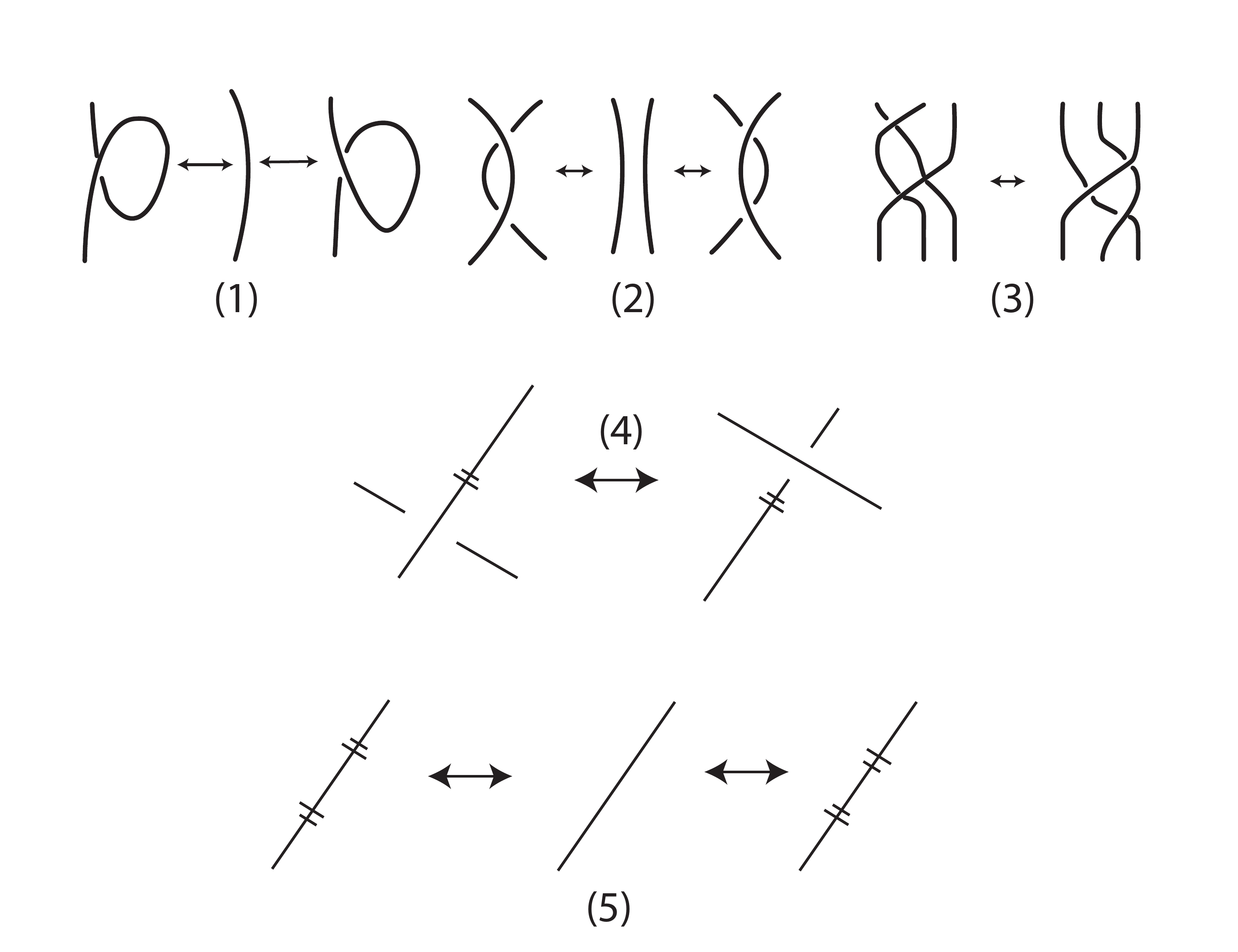}

\end{center}

 \caption{Moves for diagrams on $S_{g}$}\label{moves1}
\end{figure}
\end{prop}

\begin{cor}
Let $(L, S_{g} \times S^{1})$ and $(L', S_{g-1} \times S^{1})$ be two links in $S_{g} \times S^{1}$ and $S_{g-1} \times S^{1}$, respectively. Let $D_{L}$ and $D_{L'}$ be diagrams on $S_{g}$ and $S_{g-1}$ of $L$ and $L'$, respectively. Then $(L', S_{g-1} \times S^{1})$ is obtained from $(L, S_{g} \times S^{1})$ by a destabilization if and only if $(D_{L'}, S_{g-1})$ can be obtained from $(D_{L}, S_{g})$ by a destabilization of $S_{g}$.
\end{cor}

\begin{figure}[h]
\begin{center}
 \includegraphics[width = 12cm]{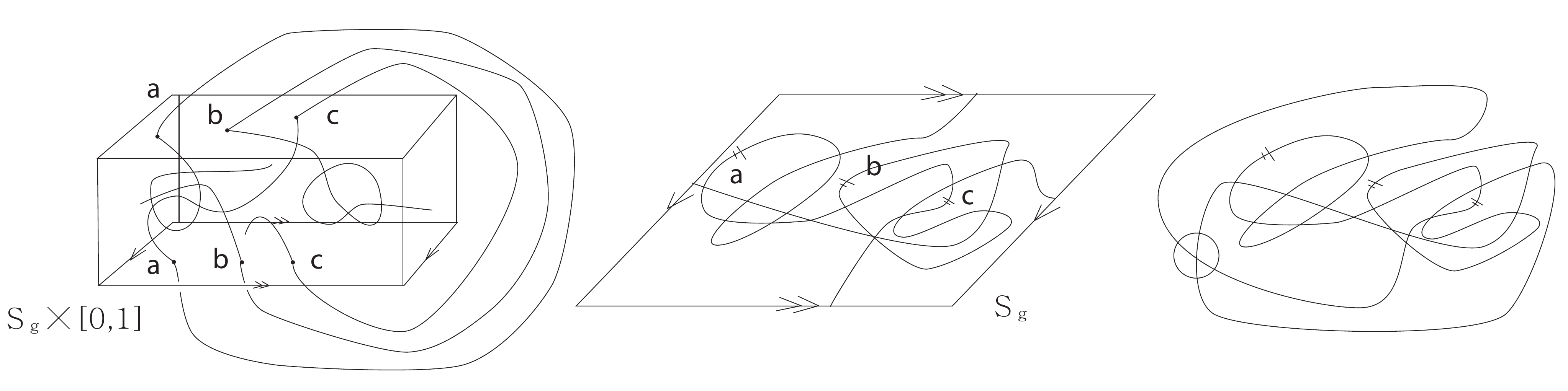}

\end{center}

\caption{Schematic figures of links in $S_{g} \times S^{1}$ and its projection on the plane}\label{fig:SS-diag-rel-2}
\end{figure}

Now, let us construct diagrams for links in $S_{g}\times S^{1}$ on the plane by using diagrams on $S_{g}$. For a link in $S_{g} \times S^{1}$ let a diagram $D$ on $S_{g}$ with double lines be given. We may assume that the diagram is drawn on $2g$-gon presentation of $S_{g}$ as illustrated in the middle of Fig.~\ref{fig:SS-diag-rel-2}. Connect points on boundaries of $2g$-gon with corresponding points by arcs outside $2g$-gon. By changing intersections between arcs outside $2g$-gon to virtual crossings, we obtain a diagram with double lines and virtual crossings, see the right in Fig.~\ref{fig:SS-diag-rel-2}. We call it {\it a diagram with double lines on the plane} for the link in $S_{g}\times S^{1}$, or simply {\it a diagram for the link in $S_{g}\times S^{1}$}.

The following theorem holds.
\begin{prop}[Kim (2018) \cite{Kim}]\label{thm:diag_on_plane}
   Let $(L, S_{g} \times S^{1})$ and $(L', S_{g'} \times S^{1})$ be two links. Let $D_{L}$ and $D_{L'}$ be diagrams of $L$ and $L'$ on the plane, respectively. Then $L$ and $L'$ are equivalent if and only if $D_{L'}$ can be obtained from $D_{L}$ by applying finitely many moves in Fig.~\ref{moves2}.

  \begin{figure}[h!]
\begin{center}
 \includegraphics[width = 12cm]{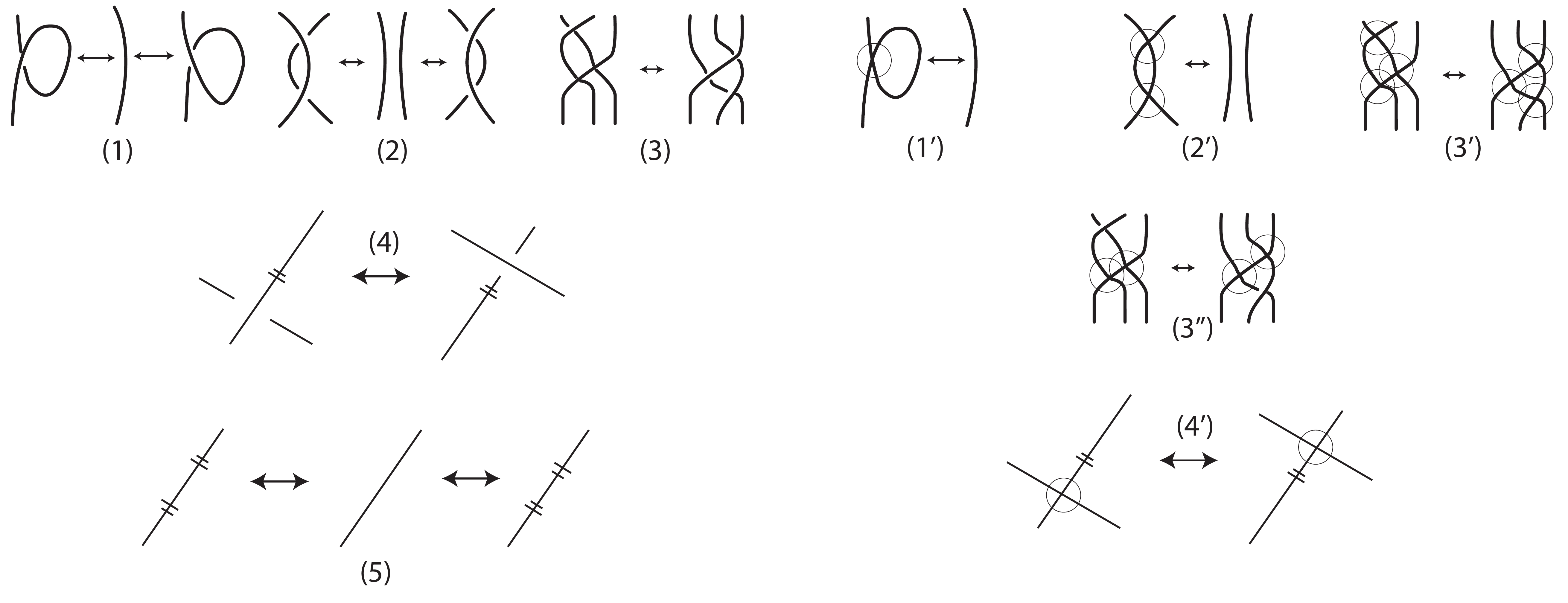}

\end{center}

 \caption{Moves for links in $S_{g}\times S^{1}$}\label{moves2}
\end{figure}
\end{prop}

{\bf Terminology.} In this paper, arcs of a diagram with double lines between two double lines are called {\em long arcs}, arcs between classical crossings are called {\em arcs} and arcs between a classical crossing and a double line are called {\em short arcs}.

\begin{rem}
As described in Fig. \ref{cro_change}, by adding two double lines one can change over/under information of a crossing. 

  \begin{figure}[h!]
\begin{center}
 \includegraphics[width = 6cm]{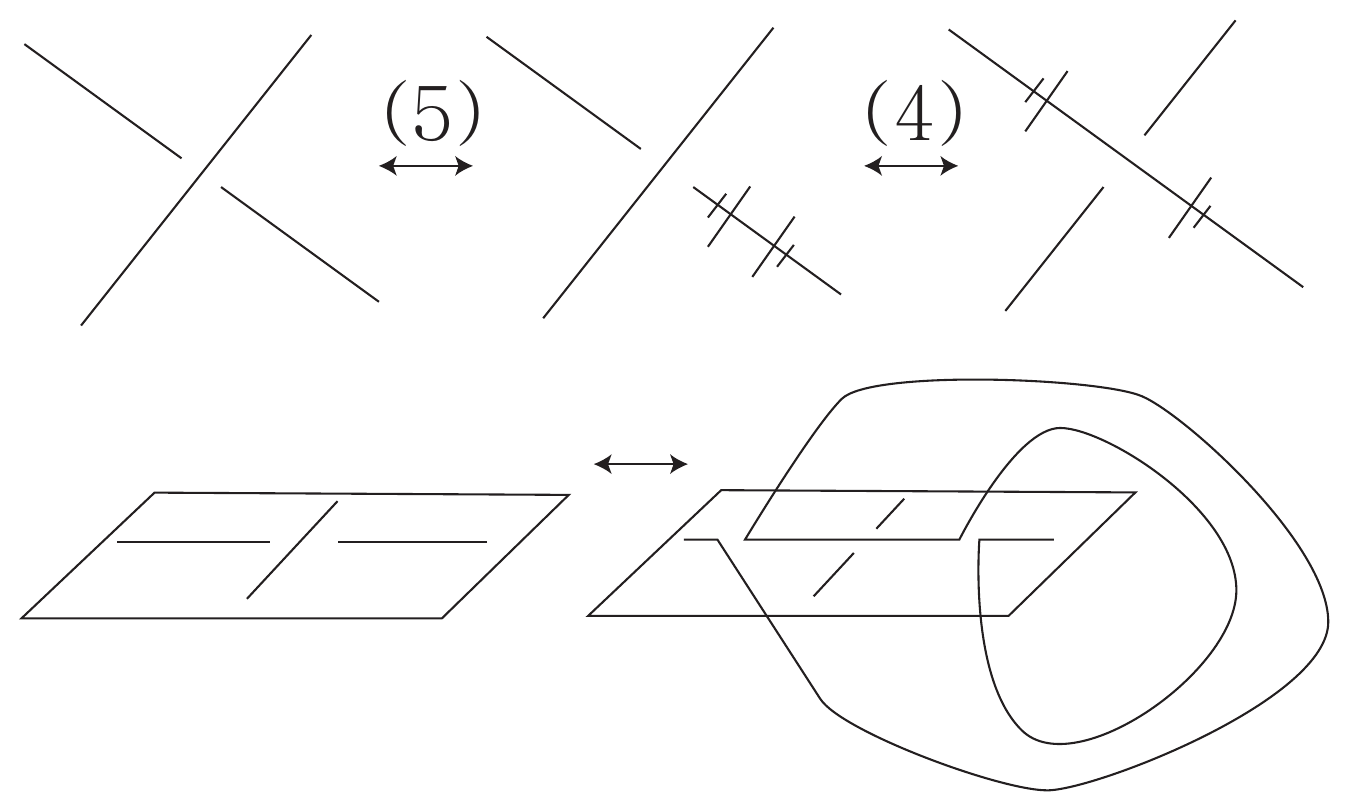}

\end{center}

 \caption{A crossing change with additional two double lines}\label{cro_change}
\end{figure}
\end{rem}


\subsection{Degree of an oriented knot in $S_{g} \times S^{1}$ and heights of long arcs}
From now on we are mainly interested in {\it oriented knots} in $S_{g} \times S^{1}$. 

Let $\Pi: \mathbb{R} \rightarrow S^{1}$ be the natural covering defined by $\Pi(r) = e^{2\pi r i}$. Then the function $Id_{S_{g}} \times \Pi : S_{g}\times \mathbb{R} \rightarrow S_{g}\times S^{1}$ is also a covering over $S_{g}\times S^{1}$, where $Id_{S_{g}} : S_{g} \rightarrow S_{g}$ is the identity map.

\begin{center}
\begin{tikzcd} 
&|[alias=Z]| S_{g}\times \mathbb{R} \arrow[r,"\phi_{2}"]\arrow[d, "Id_{S_{g}}\times \Pi"]& \mathbb{R} \arrow[d,"\Pi"]\\
S^{1}\arrow[to=Z, "\hat{K}"]\arrow[r, "K"] &S_{g}\times S^{1} \arrow{r}&S^{1}
\end{tikzcd}
\end{center}

Let $K : [0,1] \rightarrow S_{g}\times S^{1}$ be an oriented knot with $K(0)=K(1)$.
 Let $\hat{K}$ be a lifting of $K$ into $S_{g} \times \mathbb{R}$ along the covering $Id_{S_{g}} \times \Pi : S_{g}\times \mathbb{R} \rightarrow S_{g}\times S^{1}$. When $\phi_{2} \circ \hat{K}(0) =0$ for the projection $\phi_{2} : S_{g} \times \mathbb{R} \rightarrow \mathbb{R}$, {\em the degree $deg(K)$ of a knot $K$ in $S_{g} \times S^{1}$} is defined by 
 $$deg(K) := \phi_{2} \circ \hat{K}(1)\in \mathbb{Z}.$$
 It is easy to see that the degree $deg(K)$ of a knot $K$ in $S_{g} \times S^{1}$ is an invariant for knots in $S_{g} \times S^{1}$.

\begin{figure}[h!]
\begin{center}
 \includegraphics[width = 5cm]{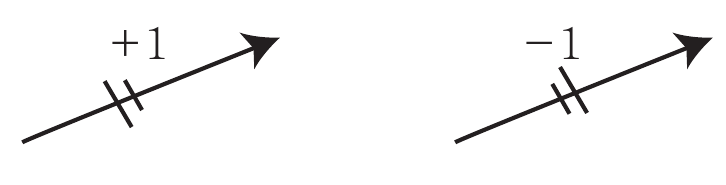}
\caption{Signs of double lines}\label{fig:signs-doublelines}
\end{center}
\end{figure}

\begin{figure}[h!]
\begin{center}
 \includegraphics[width = 8cm]{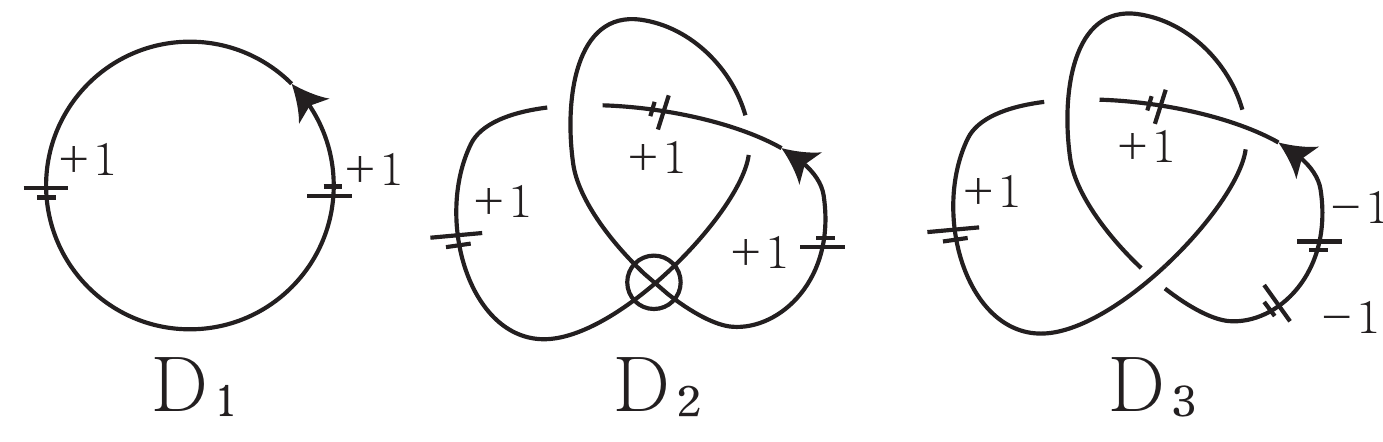}
\caption{$D_{1}$, $D_{2}$ and $D_{3}$ are oriented diagrams with double lines. Each double line has a sign. We obtain that $deg(D_{1}) = 2$, $deg(D_{2}) = 3$, $deg(D_{3}) = 0$.}\label{fig:exa-diag-degree}
\end{center}
\end{figure}

 The degree of a given oriented knot $K$ in $S_{g} \times S^{1}$ can be calculated by using a diagram with double lines in the following way:
Let $D$ be an oriented diagram with double lines of an oriented knot $K$ in $S_{g} \times S^{1}$. Let us give $\pm 1$ to double lines with respect to the orientation as described in Fig.~\ref{fig:signs-doublelines}. We call it {\em the sign of a double line}. 
Then the degree of $K$ is equal to the sum of signs of all double lines, for example, see Fig.~\ref{fig:exa-diag-degree}.

 \subsection{Heights of long arcs of $D$}
 Let $K$ be an oriented knot in $S_{g} \times S^{1}$ with degree $deg(K)=n$. Without loss of generality we may assume that the point $K(0) = K(1)$ corresponds to a double line of a diagram $D$.
 
 For a long arc $l$ of a diagram $D$ of $K$, there is a line segment $l'$ in the lifting $\hat{K}$ to $S_{g} \times \mathbb{R}$ with $\hat{K}(0) \in S_{g} \times \{0\}$ corresponding to $l$. If $\phi_{2}(l') \in (a, a+1)$ for some $a \in \mathbb{Z}$, then we give a label $a \in \mathbb{Z}_{n}$ (or $a\in \mathbb{Z}$ when $n=0$) to $l$ and we call it {\em the height of the long arc $l$}.

\begin{figure}[h!]
\begin{center}
 \includegraphics[width = 6cm]{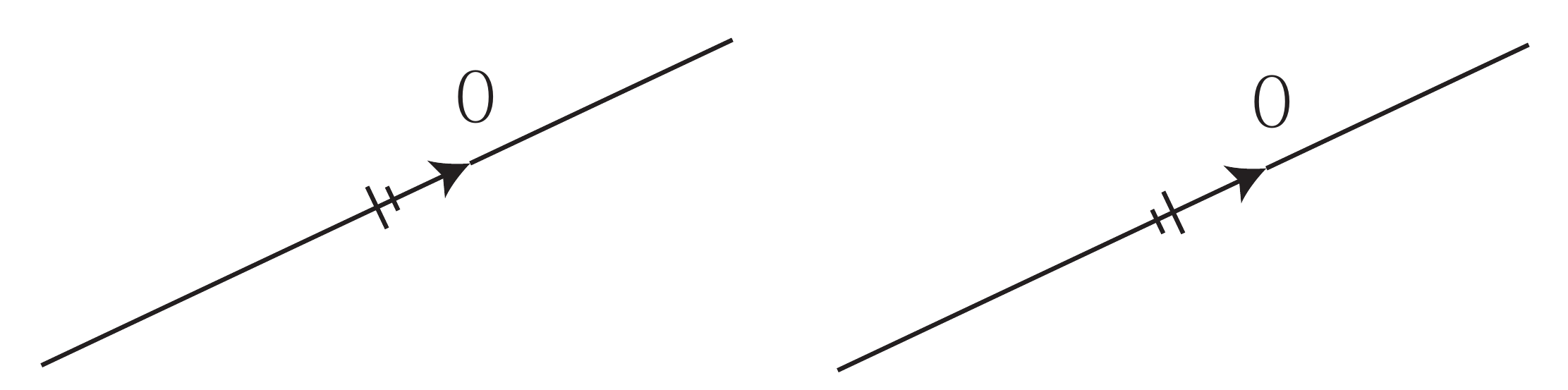}

\end{center}
 \caption{Starting point determining heights of long arcs}\label{fig:rule-height1}
\end{figure}

\begin{exa}
Let $K$ be an oriented knot in $S^{2} \times S^{1}$ as described in Fig.~\ref{exa_label_arc}, where $S^{2}$ is a 2-dimensional sphere. The knot $K$ has the degree $2$ and it has a diagram $D_{K}$ of the trivial knot diagram with two double lines. One can see that the arc of $K$ colored by red corresponds to the arc of $\hat{K}$ placed in $S^{2} \times (0,1)$, but the arc of $K$ colored by green corresponds to the arc of $\hat{K}$ placed in $S^{2} \times (1,2)$. Note that the red and green arcs of $K$ correspond to the long arcs of $D_{K}$ colored by red and green respectively. Now we give numbers $0$ and $1$ to red and green long arcs of $D_{K}$ respectively. 
  \begin{figure}[h!]
\begin{center}
 \includegraphics[width =8cm]{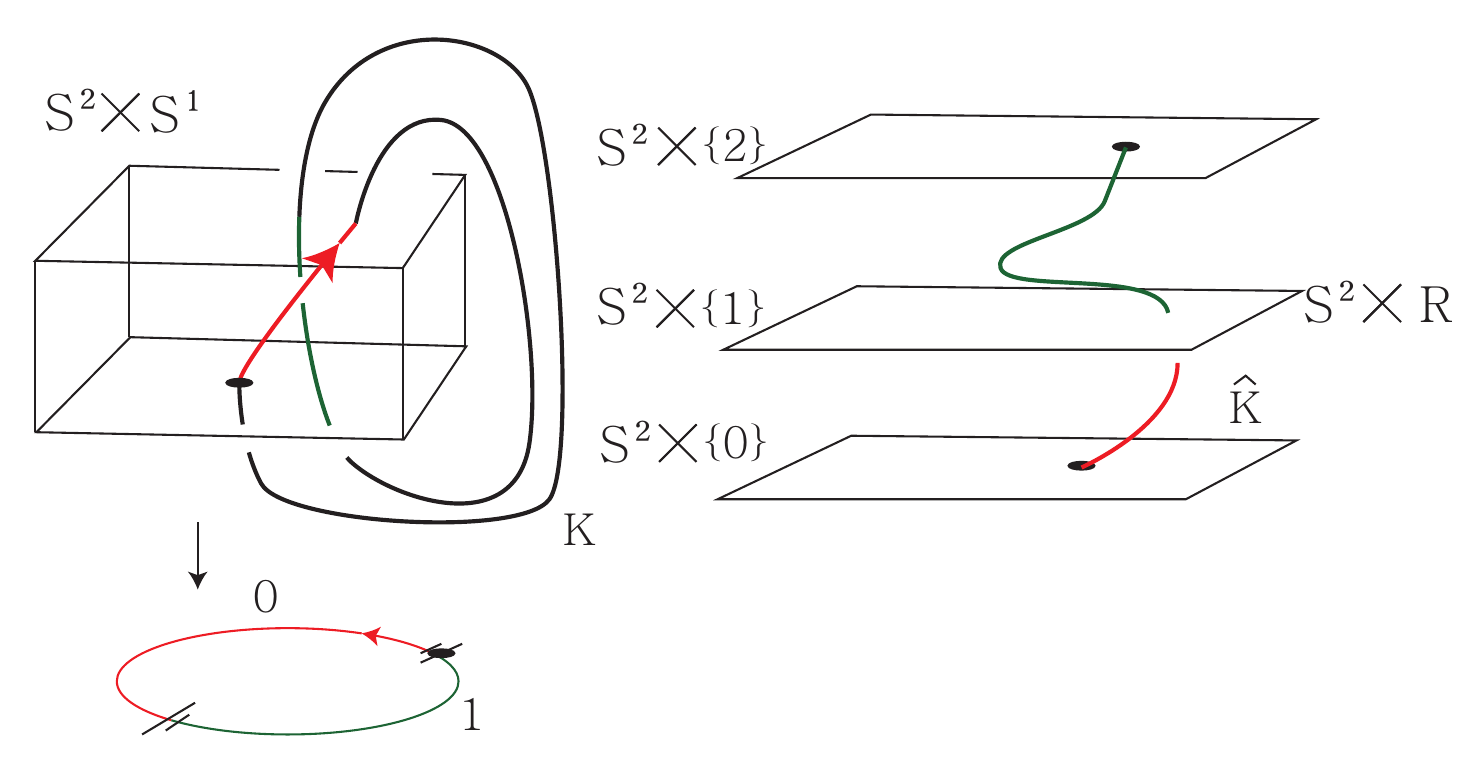}

\end{center}
 \caption{A knot in $S^{2} \times S^{1}$ with degree $2$}\label{exa_label_arc}
\end{figure}
\end{exa}
 Heights of long arcs of a diagram $D$ of $K$ can be determined in combinatoric way: Let $K$ be an oriented knot in $S_{g} \times S^{1}$ with the degree $deg(K) = n$. Let $D$ be an oriented diagram with double lines of $K$. Let us give signs for double lines. Take a double line as a starting point and give a height $0$ to the long arc which is the following long arc with respect to the orientation of $D$, see Fig.~\ref{fig:rule-height1}.

\begin{figure}[h!]
\begin{center}
 \includegraphics[width = 6cm]{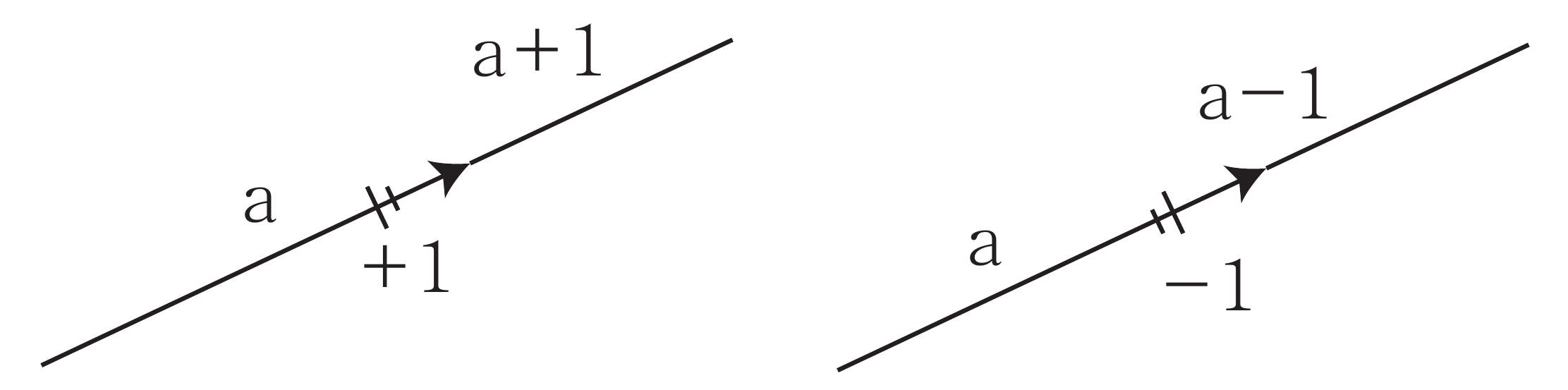}

\end{center}
 \caption{Heights of two long arcs with common double line}\label{fig:rule-height2}
\end{figure}

Following the diagram respecting the given orientation when we meet a double line with height $a \in \mathbb{Z}_{n}$ coming from a long arc, we give $a+\epsilon \in \mathbb{Z}_{n}$ to the next long arc, where $\epsilon$ is the sign of the double line, see Fig.~\ref{fig:rule-height2}.

\begin{exa}
Let $D$ be an oriented diagram with double lines in Fig.~\ref{fig:exa_height_diag_deg0} of an oriented knot $K$ in $S_{g} \times S^{1}$. The degree of $K$ in Fig.~\ref{fig:exa_height_diag_deg0} is $0$. With the fixed starting point at a double line, we obtain heights of long arcs as the right figure in Fig.~\ref{fig:exa_height_diag_deg0}.
\begin{figure}[h!]
\begin{center}
 \includegraphics[width = 10cm]{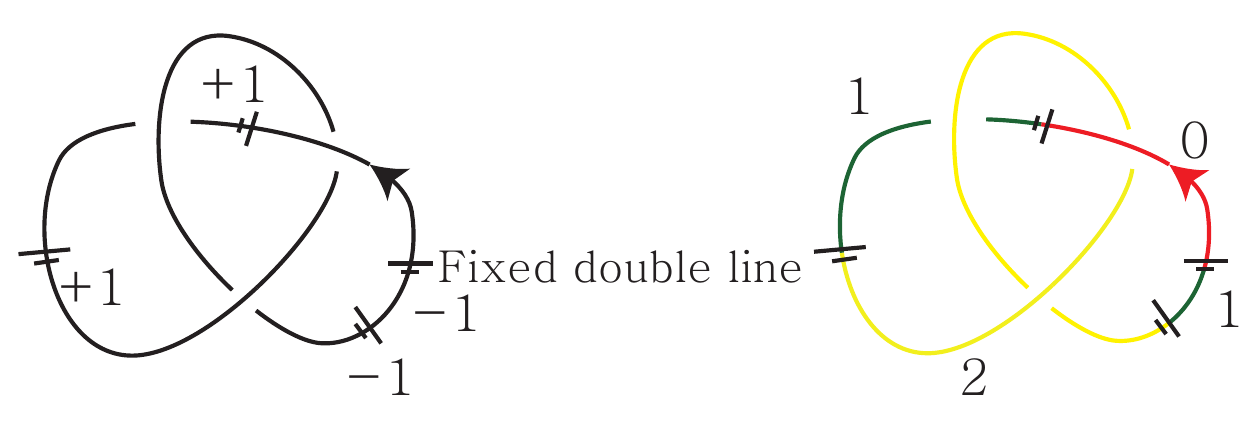}
\caption{Heights  on long arcs of a knot $K$ of degree $0$ valued in $\mathbb{Z}_{0}=\mathbb{Z}$}\label{fig:exa_height_diag_deg0}
\end{center}
\end{figure}

\end{exa}

\begin{exa}
Let $D$ be an oriented diagram with double lines in Fig.~\ref{fig:exa_height_diag_deg3} of an oriented knot in $S_{g} \times S^{1}$. The degree of $K$ is $3$. With the fixed starting point at a double line, we obtain heights of long arcs. Note that the long arc with height $0$ has height $3$ simultaneously by the rule of determining heights in combinatoric way. Since heights are valued in $\mathbb{Z}_{3}$, it does not make any ambiguity.

\begin{figure}[h!]
\begin{center}
 \includegraphics[width = 8cm]{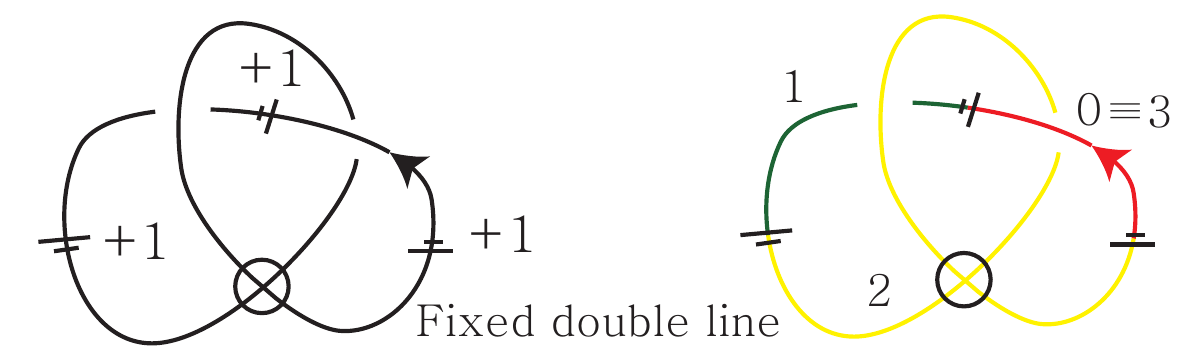}
\caption{Heights on long arcs of a knot $K$ of degree $3$ valued in $\mathbb{Z}_{3}$}\label{fig:exa_height_diag_deg3}
\end{center}
\end{figure}

\end{exa}

\begin{rem}
 Let $\mathcal{A}_{D}$ be the set of long arcs of $D$. Then heights of long arcs can be presented by a map from $\mathcal{A}_{D}$ to $\mathbb{Z}_{n}$. For each choice of a double line $d$ as a starting point we obtain $\phi_{d}: \mathcal{A}_{D}\rightarrow \mathbb{Z}_{n}$. For a pair of double lines $d_{1}$ and $d_{2}$, there exists a map $t_{s} : \mathbb{Z}_{n} \rightarrow \mathbb{Z}_{n}$ defined by $t_{s}(x) = x+s$ such that $\phi_{d_{2}}= t_{s} \circ \phi_{d_{1}}$. \label{rem:height-translation}
\end{rem}
\section{Liftings of knots in $S_{g} \times S^{1}$}
In the previous section we define the degree of an oriented knot in $S_{g} \times S^{1}$ by using the lifting of the knot along $Id_{S_{g}} \times \Pi : S_{g}\times \mathbb{R} \rightarrow S_{g}\times S^{1}$. Since a knot $K$ is an embedding of $S^{1}$ into $S_{g}\times S^{1}$, there exists a lifting of $K$ to $S_{g} \times \mathbb{R}$ along the projection $Id_{S_{g}} \times \Pi$. Since $\mathbb{R} \cong (0,1) \subset [0,1]$, the lifting can be considered as an arc in thickened surface $S_{g} \times [0,1]$ and one can expect that it can be presented by a diagram on the plane. In the present section we discuss how to obtain the diagram of the liftings step by step by using diagrams with heights of long arcs.

\subsection{Lifting of knots in $S_{g} \times S^{1}$ with degree $0$}
Let $K$ be an oriented knot in $S_{g} \times S^{1}$ with degree $0$. Let $\hat{K}$ be a lifting of $K$ to $S_{g} \times \mathbb{R}$. Since $K$ has the degree $0$, $\hat{K}$ is a knot in $S_{g} \times \mathbb{R}$, that is, $\hat{K}$ is a virtual knot. 
To obtain a diagram of $\hat{K}$ let us visualize the lifting $\hat{K}$ to $S_{g} \times \mathbb{R}$. It consists of 3 steps:\\

\textbf{Step 1.} Let $D$ be an oriented diagram with double lines of $K$. Let us fix a double line on a diagram and give a height for each long arc as done in the previous section. Let us say that the smallest height is $m$ and the biggest height is $M$.\\

\begin{figure}[h!]
\begin{center}
 \includegraphics[width = 6cm]{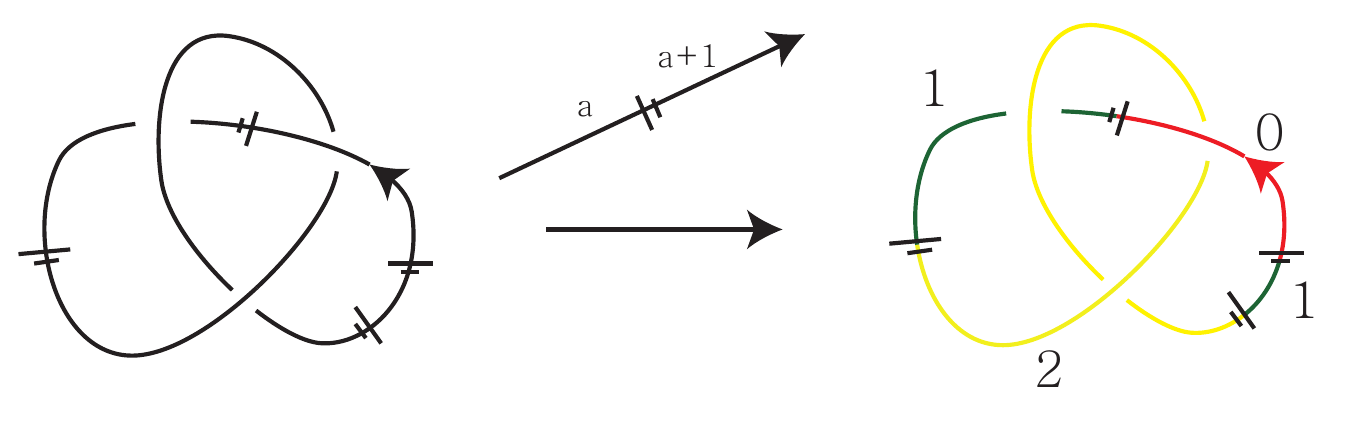}
 \caption{Step 1 for knots of degree $0$}\label{fig:exa_label_diag_deg0}
\end{center}
\end{figure}

\begin{figure}[h!]
\begin{center}
 \includegraphics[width = 6cm]{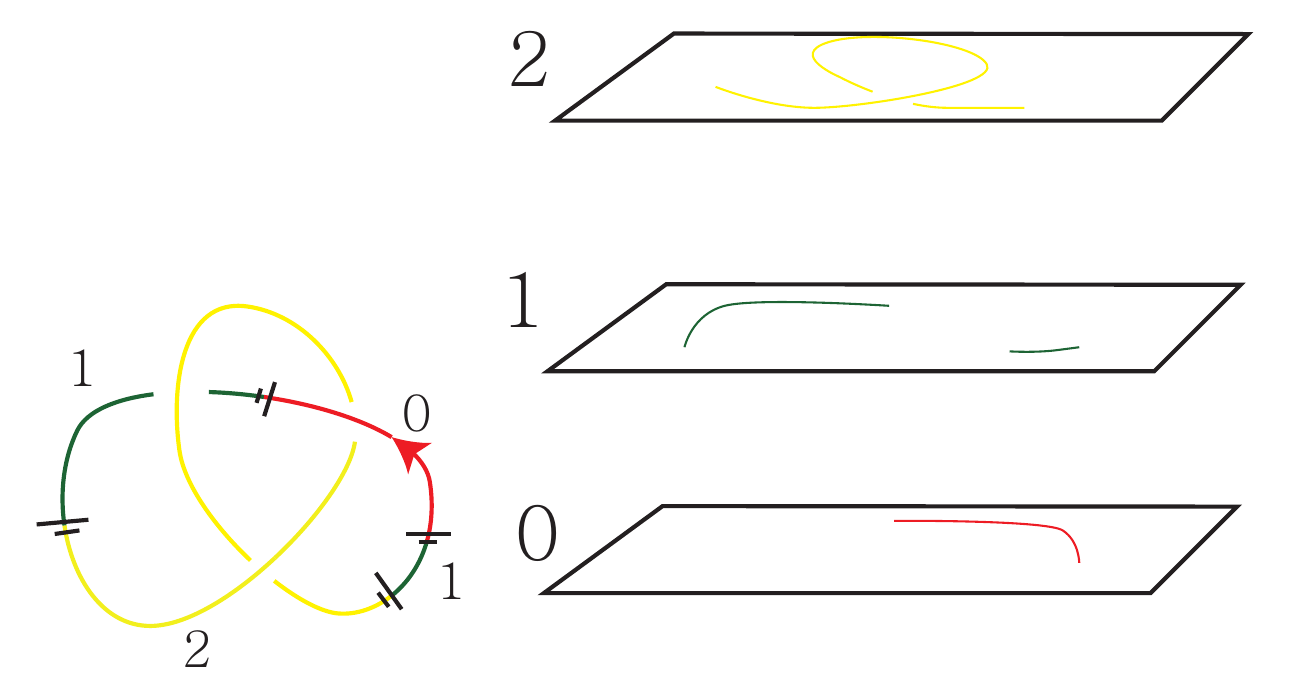}
\end{center}
\caption{Step 2 for knots of degree $0$}\label{fig:algo_lifting_deg0-1}
\end{figure}

\textbf{Step 2.} Let us imagine $M-m+1$ parallel planes placed as described in Fig.~\ref{fig:algo_lifting_deg0-1}. Give numbers from bottom to top of planes by integers from $m$ to $M$. Draw $M-m+1$ copies of $D$ for each plane. For a copy of a diagram on each plane with a number $s$ erase long arcs which have not the height $s$. \\

\textbf{Step 3.} 
Let us suppose that we are walking on a long arc of $D$ with height $s$. When we meet a double line, if it is a longer line, then we connect the corresponding arc on the plane with number $s$ to the arc on the plane with number $s+1$, which is the arc corresponding to the long arc adjacent to the double line, which we met. If it is a shorter line, then we connect it to the arc on the plane with number $s-1$ corresponding to the long arc adjacent to the double line, which we met. Since $deg(K)=0$, we must come back to the starting point and it is a virtual knot in $S_{g} \times \mathbb{R}$. 

\begin{figure}[h!]
\begin{center}
 \includegraphics[width = 8cm]{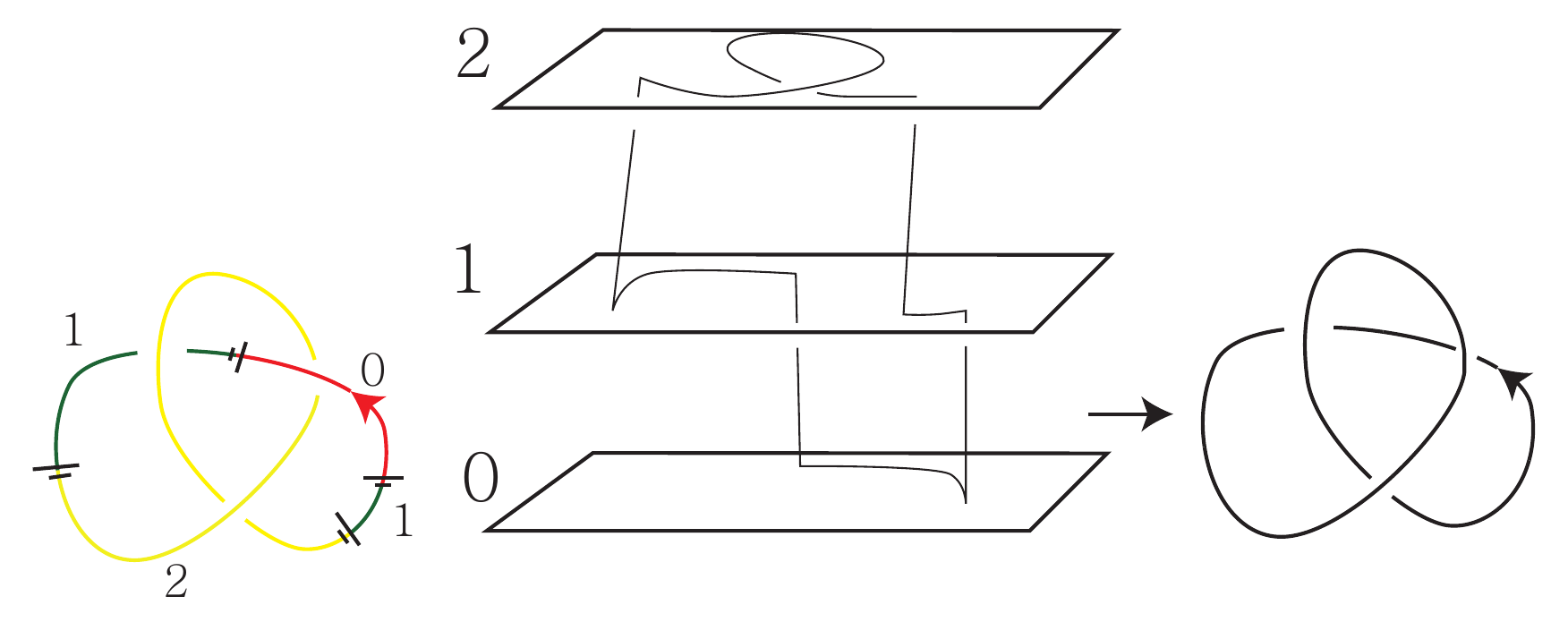}
\end{center}
\caption{Step 3 for knots of degree $0$}\label{fig:algo_lifting_deg0-3}
\end{figure}

From our visualization of $\hat{K}$ one can obtain a virtual diagram $\hat{D}$ as follows: for a diagram $D$ with double lines, we give heights to long arcs of $D$. For each crossing, if two heights of over/under crossings are same, then we remain it. If two heights of over/under crossings are different, then we change (if it is necessary) over/under information so that the arc with bigger height becomes over crossing. By removing double lines, we obtain a virtual knot diagram $\hat{D}$ corresponding to the lifting $\hat{K}$.


\begin{lem}\label{lem:lifting-deg0}
Let $D$ and $D'$ be two oriented diagrams with double lines of knots in $S_{g} \times S^{1}$. If they are equivalent, then $\hat{D}$ and $\hat{D'}$ are equivalent as virtual knots.
\end{lem}

\begin{proof}

Assume that $D'$ is obtained from $D$ by applying one of the moves in Fig.~\ref{moves2}. If $D'$ is obtained from $D$ by applying moves (1), (2), (3), (1'), (2'), (3') and (3''), then it is easy to see that $\hat{D'}$ can be obtained from $\hat{D}$ by applying virtual Reidemeister moves.

Suppose that $D'$ is obtained from $D$ by applying the move (4). The line segments in move (4) have heights $a,a+1$ and $b$ as described in Fig.~\ref{lifting4_type1} and Fig.~\ref{lifting4_type2}. As shown in Fig.s, when we lift diagrams, the difference is just a change of the place of a classical crossing and it follows that $\hat{D} \cong \hat{D}'$. 

  \begin{figure}[h!]
\begin{center}
 \includegraphics[width = 12cm]{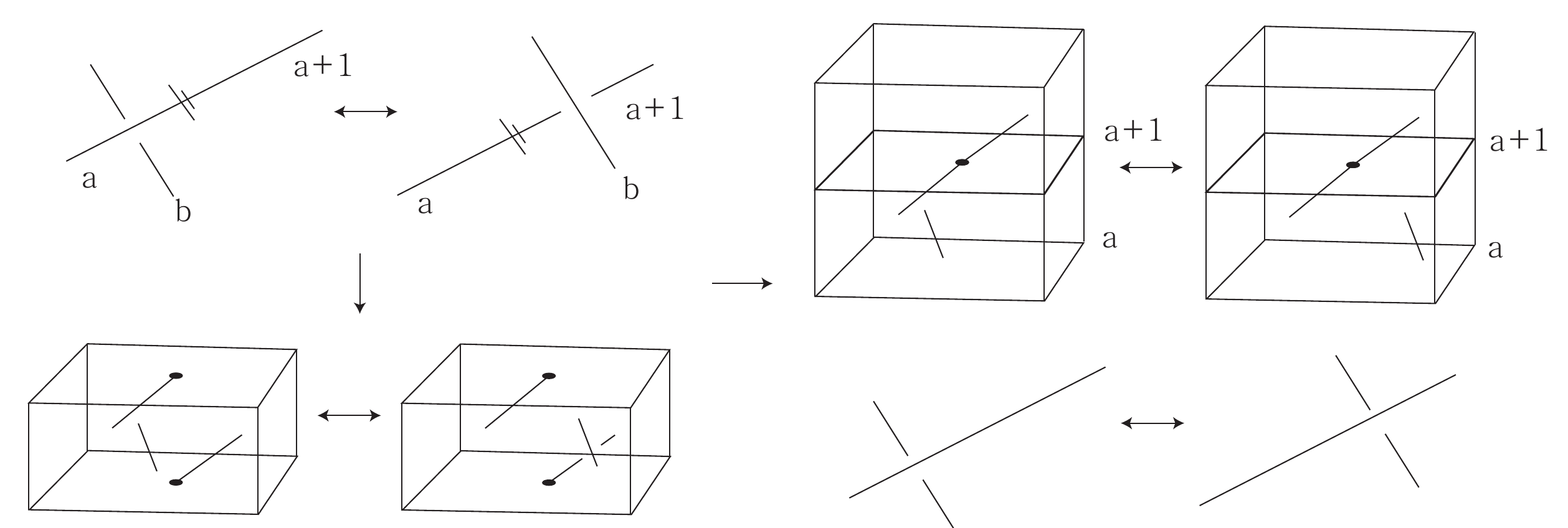}

\end{center}
 \caption{ $a\geq b$}\label{lifting4_type1}
\end{figure}

  \begin{figure}[h!]
\begin{center}
 \includegraphics[width = 12cm]{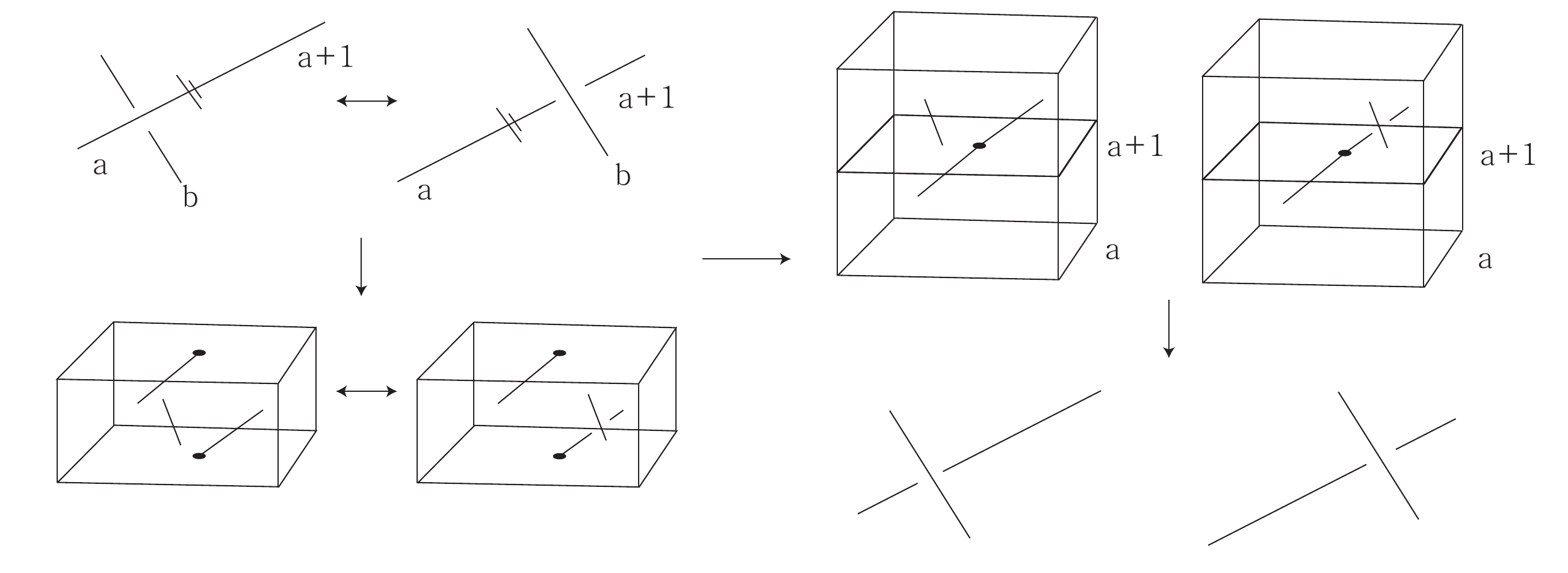}

\end{center}
 \caption{ $a< b$}\label{lifting4_type2}
\end{figure}

Suppose that $D'$ is obtained from $D$ by applying the move (5), see Fig.~\ref{lifting5}. When we lift the link to $S_{g} \times \mathbb{R}$ the line segment between two double-lines and others are placed in the different levels. Note that in $S_{g} \times \mathbb{R}$ under the line segment corresponding to the long arc between two double-lines there are no other arcs, so we can push it down as described under of Fig.~\ref{lifting5}. That is, $\hat{D} \cong \hat{D}'$ and the proof is completed.

  \begin{figure}[h!]
\begin{center}
 \includegraphics[width = 12cm]{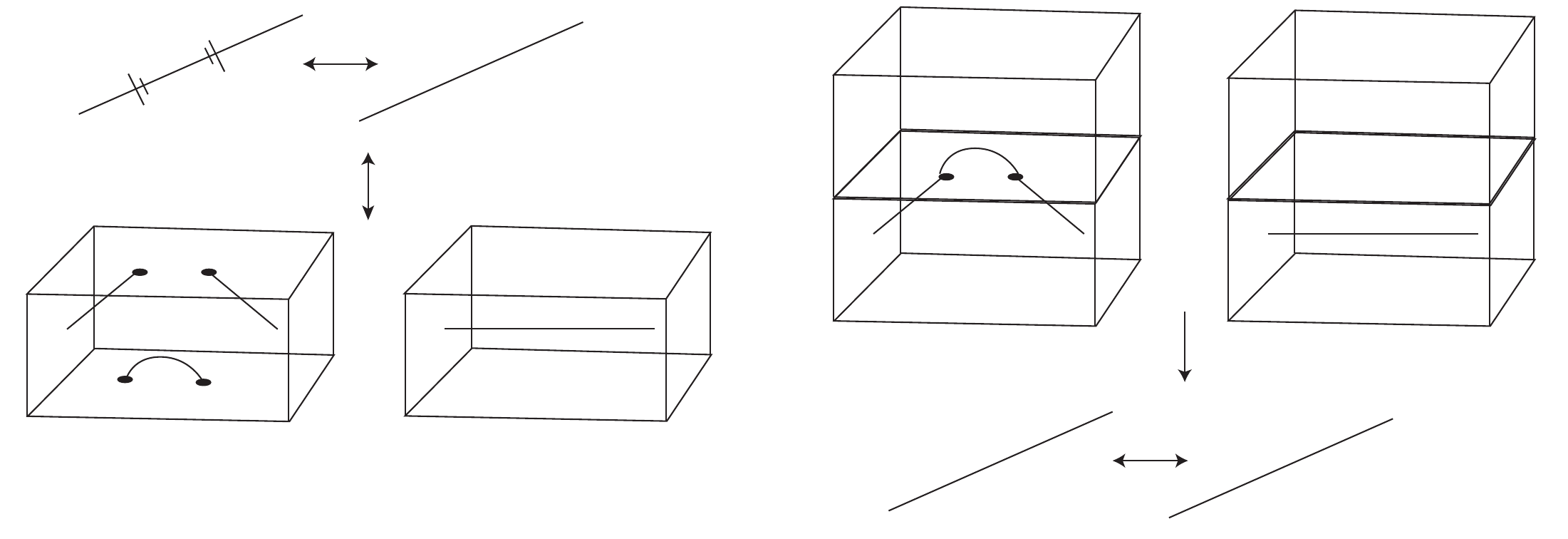}

\end{center}
 \caption{Move (5) and the corresponding liftings}\label{lifting5}
\end{figure}

\end{proof}

\begin{rem}
If we consider liftings $\hat{K}_{s}$ such that $\hat{K}_{s}(0) = \hat{K}_{s}(0) = (x, s)$, then we obtain a link of infinitely many components. Note that the classical crossing between two long arcs with heights $a$ and $b$, then it corresponds to the crossing of the link $\sqcup_{s \in \mathbb{Z}} \hat{K}_{s}$ of infinite components between two components $\hat{K}_{s}$ and $\hat{K}_{s+b-a}$, see Fig. \ref{exa_inf-lifting_deg0}.
 \begin{figure}[h!]
\begin{center}
\includegraphics[width = 8cm]{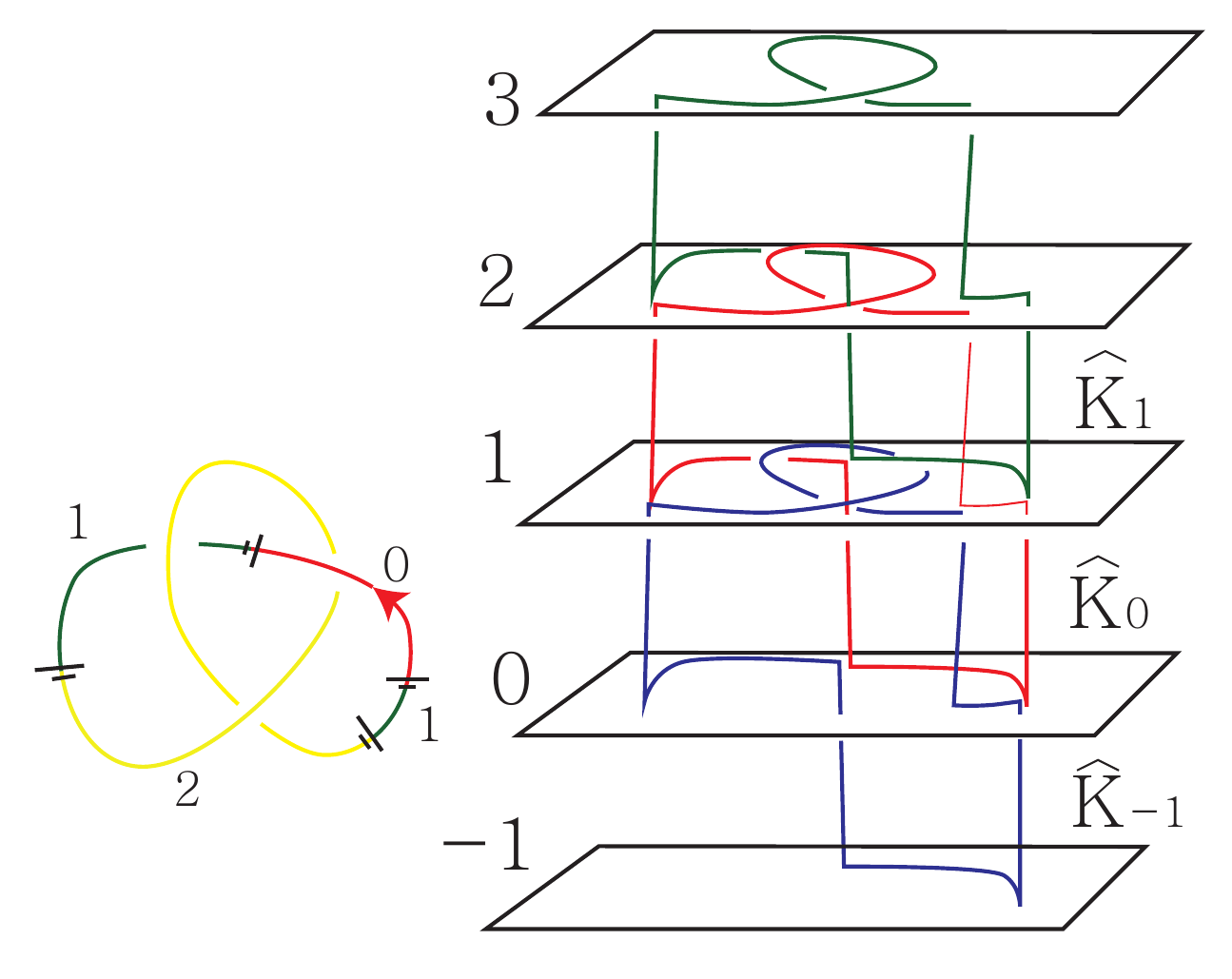}
\end{center}
\caption{Liftings $\{\hat{K}_{s}\}_{s \in \mathbb{Z}}$ of a knot $K$ of degree $0$ where $\hat{K}_{s}(0) = (x,s)$}\label{exa_inf-lifting_deg0}
\end{figure}
Moreover, if we restrict the number of components of the lifting, then we obtain a link in $S_{g} \times \mathbb{R}$ of $m$ components as a lifting of a knot in $S_{g} \times S^{1}$.
\end{rem}

\begin{cor}\label{cor:multi-lifting}
Let $K$ and $K'$ be knots in $S_{g} \times S^{1}$. If $K$ and $K'$ are equivalent as knots in $S_{g} \times S^{1}$, then the liftings $\sqcup_{s =1}^{m} \hat{K}_{s}$ and $\sqcup_{s =1}^{m} \hat{K'}_{s}$ are equivalent as virtual links.  
\end{cor}

 {\bf The construction of a diagram of $\sqcup_{s =1}^{m} \hat{K}_{s}$ of $m$ components for a knot $K$ in $S_{g} \times S^{1}$ with degree $0$} consists of 5 steps.

{\bf Step 1.} For a diagram $D$ with double lines, we draw $m$ parallels. But we do not determine over/under information yet.

\begin{figure}[h!]
\begin{center}
 \includegraphics[width = 8cm]{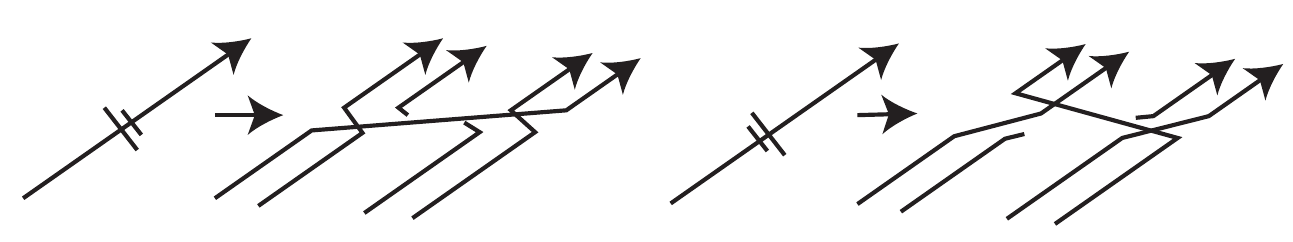}
\end{center}
 \caption{Construction of diagram}\label{diag-const-lifting-doubleline}
\end{figure}

{\bf Step 2.} Twist $m$ parallel arcs by $\pi/m$ on the part corresponding to double lines as described in Fig.~\ref{diag-const-lifting-doubleline}. Note that over/under information is not determined yet.

{\bf Step 3.} Fix a point for each component so that $m$ points are placed in one line vertical to all arcs. And give $0,\dots, m-1$ for arcs with the fixed point from left to right with respect to the given orientation. 

\begin{figure}[h!]
\begin{center}
 \includegraphics[width = 8cm]{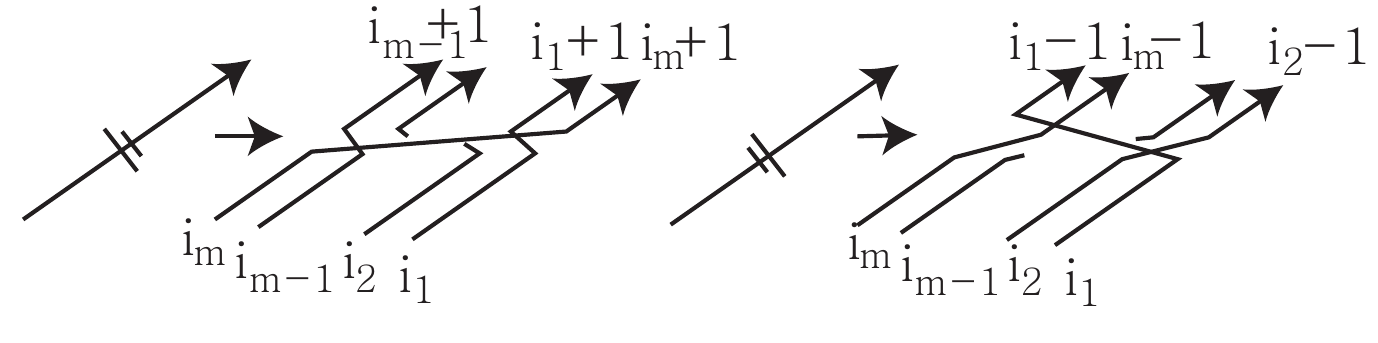}

\end{center}
 \caption{Construction of diagram}\label{diag-const-lifting-numbering}
\end{figure}

{\bf Step 4.} Walking along the diagram, we number arcs so that the number of arc is not changed when we pass through the crossings corresponding to a classical crossing of a knot $K$. When we pass classical crossings associated to double lines, the numbering of arcs before twisting is changed by $\pm 1$ with the rule described in Fig.~\ref{diag-const-lifting-numbering}, but the numbering of arcs with crossings is same to the numbering before twisting. 

\begin{figure}[h!]
\begin{center}
 \includegraphics[width = 12cm]{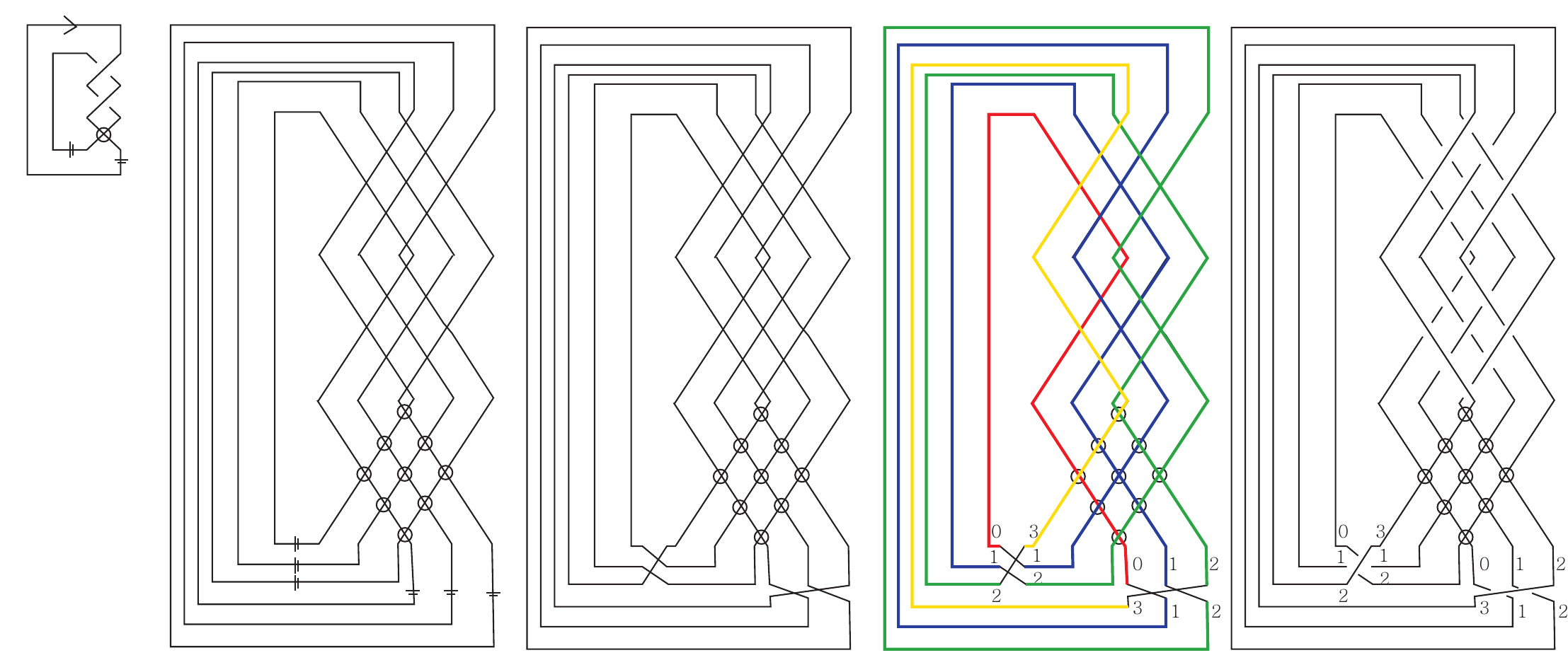}

\end{center}

\caption{Construction of the diagram of 3-fold lifting}\label{exa-construction-lift-dl}
\end{figure}

{\bf Step 5.} We change crossings associated to a virtual crossing of $D$ to virtual crossings.
For each crossing associated to a classical crossing or a double line of $D$ if arcs of a crossing have different number, then the arc with higher numbering becomes over crossing, otherwise, follow the over/under information of original diagram. And we obtain a virtual knot diagram $\sqcup_{s=1}^{m} \hat{D}_{s}$ of $m$ components.
For example, see Fig.~\ref{exa-construction-lift-dl}.

\begin{proof}[Proof of Corollary \ref{cor:multi-lifting}]
Corollary \ref{cor:multi-lifting} follows from the definition of liftings. But this can be proved by using diagrams constructed as above. Let $D$ and $D'$ be diagrams with double lines associated knots $K$ and $K'$ in $S_{g} \times S^{1}$. If $D'$ is obtained from $D$ by applying moves (1), (2), (3), (1'), (2'), (3') and (3''), then it is easy to see that $\sqcup_{s=1}^{m} \hat{D'}_{s}$ can be obtained from $\sqcup_{s=1}^{m}\hat{D}_{s}$ by applying virtual Reidemeister moves.

Suppose that $D'$ is obtained from $D$ by the move (4). Then $\sqcup_{s=1}^{m}\hat{D}_{s}$ and $\sqcup_{s=1}^{m}\hat{D'}_{s}$ have diagrams as described in Fig.~\ref{fig:multi-lifting-deg0-move4}. 
\begin{figure}[h!]
\begin{center}
\includegraphics[width = 8cm]{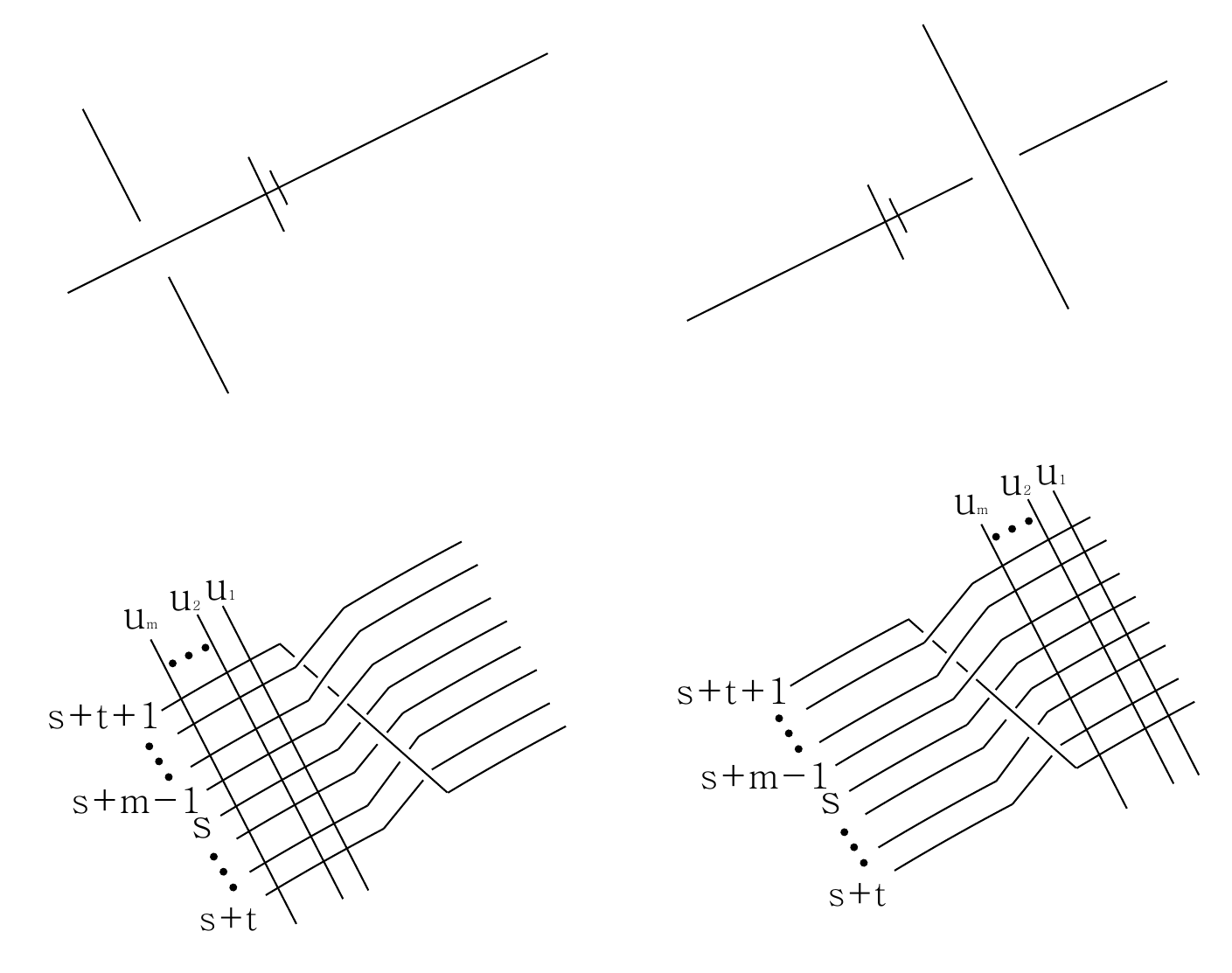}
\end{center}
\caption{Diagrams $\sqcup_{s=1}^{m}\hat{D}_{s}$ and $\sqcup_{s=1}^{m}\hat{D'}_{s}$ obtained from $D$ and $D'$ when $D'$ is obtained from $D$ by the move (4)}\label{fig:multi-lifting-deg0-move4}
\end{figure}
It is enough to show that one arc with number $u_{i}$ can pass crossings. If $s\leq u \leq s+t$, then the arc with number $u$ can pass through the crossings by Reidemeister moves as described in Fig. \ref{fig:multi-lifting-deg0-move4-case1}.
\begin{figure}[h!]
\begin{center}
\includegraphics[width = 10cm]{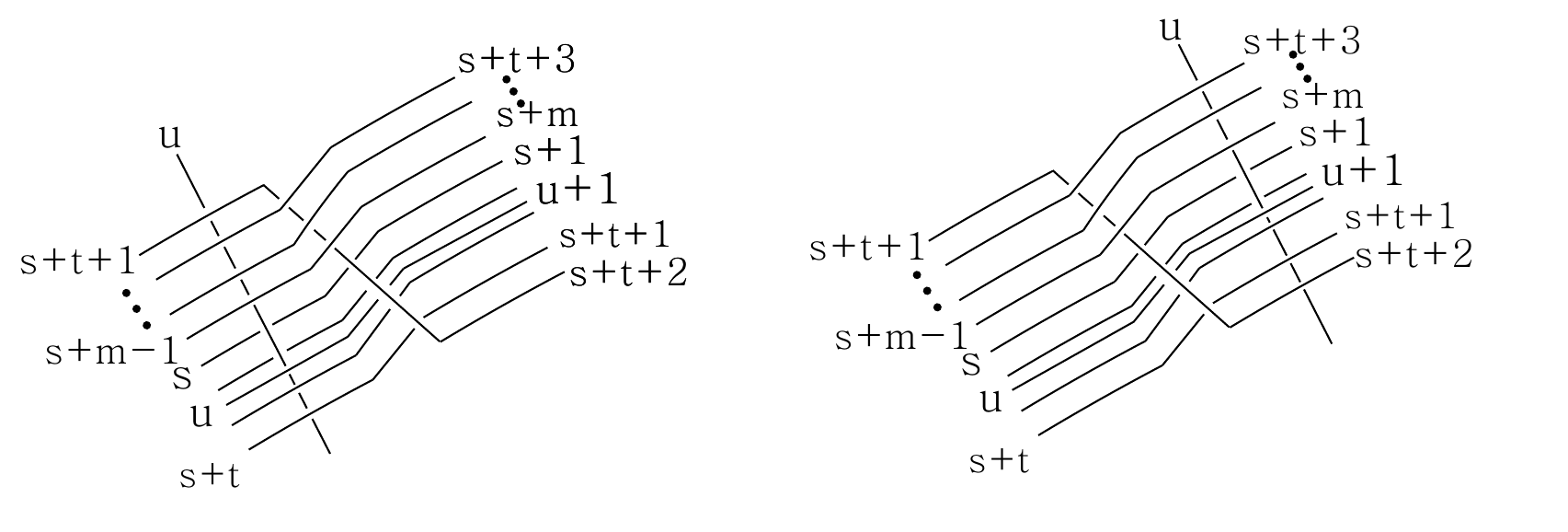}
\end{center}
\caption{Case 1. $s\leq u \leq s+t$}\label{fig:multi-lifting-deg0-move4-case1}
\end{figure}
If $s+t+1\leq u \leq m-1$, then the arc with number $u$ can pass through the crossings by Reidemeister moves as described in Fig. \ref{fig:multi-lifting-deg0-move4-case2}.
\begin{figure}[h!]
\begin{center}
\includegraphics[width = 10cm]{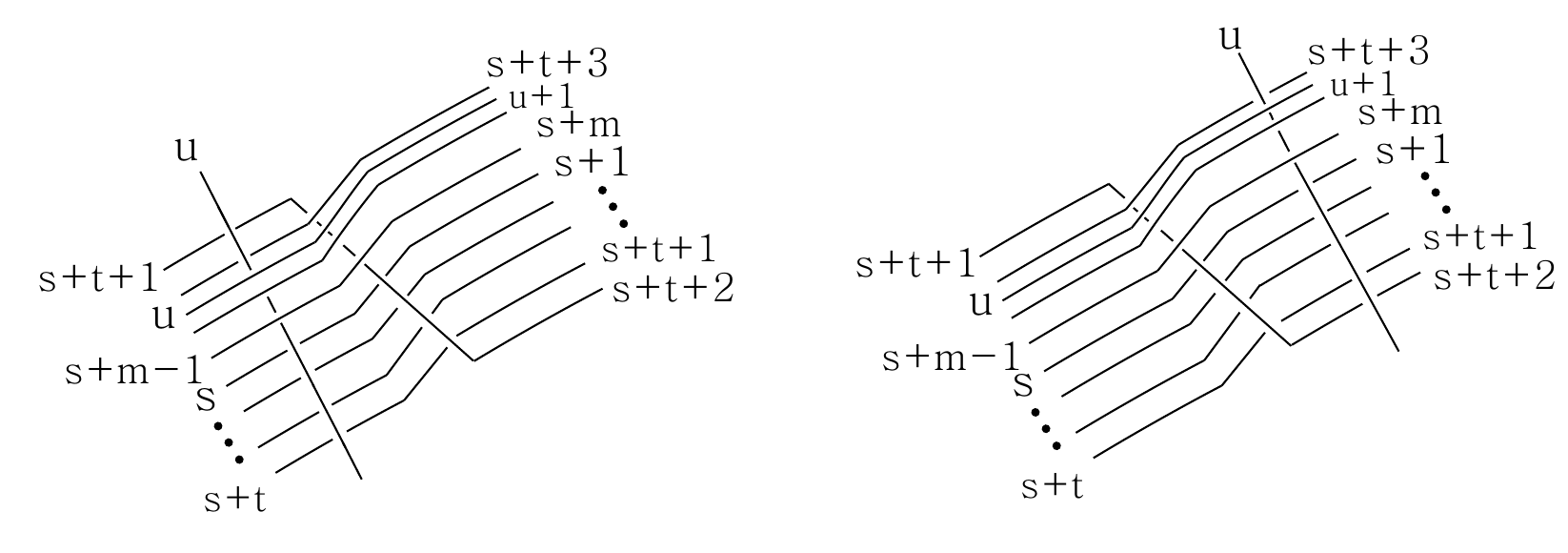}
\end{center}
\caption{Case 2. $s+t+1\leq u \leq m-1$}\label{fig:multi-lifting-deg0-move4-case2}
\end{figure}

If $D'$ is obtained from $D$ by the move (5), then $\sqcup_{s=1}^{m}\hat{D'}_{s}$ can be obtained from $\sqcup_{s=1}^{m}\hat{D}_{s}$ by a series of the second Reidemeister moves. Let us leave the details for the move (5) to readers.

\end{proof}
\subsection{Lifting of knots in $S_{g} \times S^{1}$ with degree $k\neq 0$}

Let $K$ be an oriented knot in $S_{g} \times S^{1}$ with degree $k\neq 0$. Then there exists a lifting $\hat{K}$ to $S_{g} \times \mathbb{R}$, but, unlike the case of knots in $S_{g} \times S^{1}$ with degree $0$, it is a long knot obtained by connected sum of infinitely many copies of $\hat{K} \cap [0,k]$. But, if we consider a covering $p_{k}$ from $S_{g} \times S^{1}$ to $S_{g} \times S^{1}$ defined by $p_{k}(x,z) = (x,z^{k})$, the lifting $\hat{K}$ to $S_{g} \times S^{1}$ becomes a knot in $S_{g} \times S^{1}$ of degree $1$. The algorithm to obtain a diagram of $\hat{K}$ of a knot $K$ in $S_{g}\times S^{1}$ with degree $k\neq 0$ is similar to the algorithm for the case of knots of degree $0$; it consists of 3 steps:

\textbf{Step 1.} Let $D$ be an oriented diagram with double lines of $K$ in $S_{g} \times S^{1}$ of degree $k \neq 0$. Let us fix a double line on a diagram and give a height for each long arc. We consider heights as elements in $\mathbb{Z}_{k}$.\\

\begin{figure}[h!]
\begin{center}
 \includegraphics[width = 6cm]{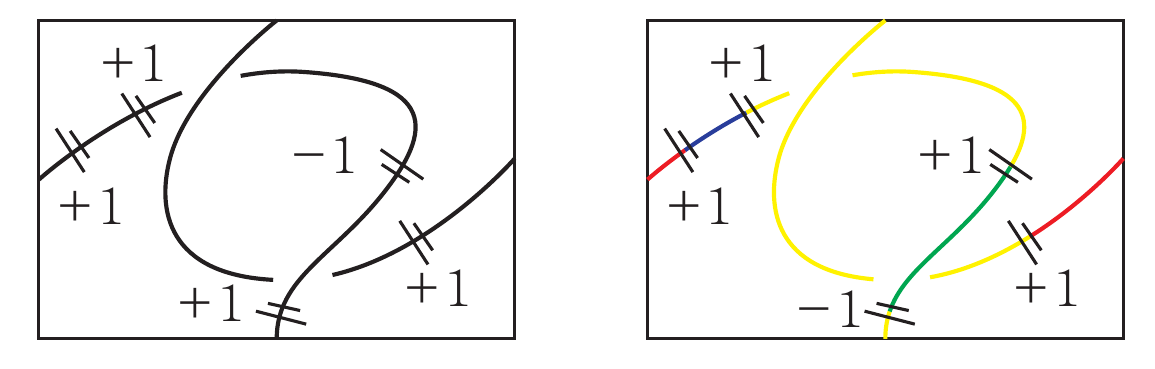}
 \caption{Step 1 for knots of degree $3$}\label{fig:exa_label_diag_deg3}
\end{center}
\end{figure}

\begin{figure}[h!]
\begin{center}
 \includegraphics[width = 6cm]{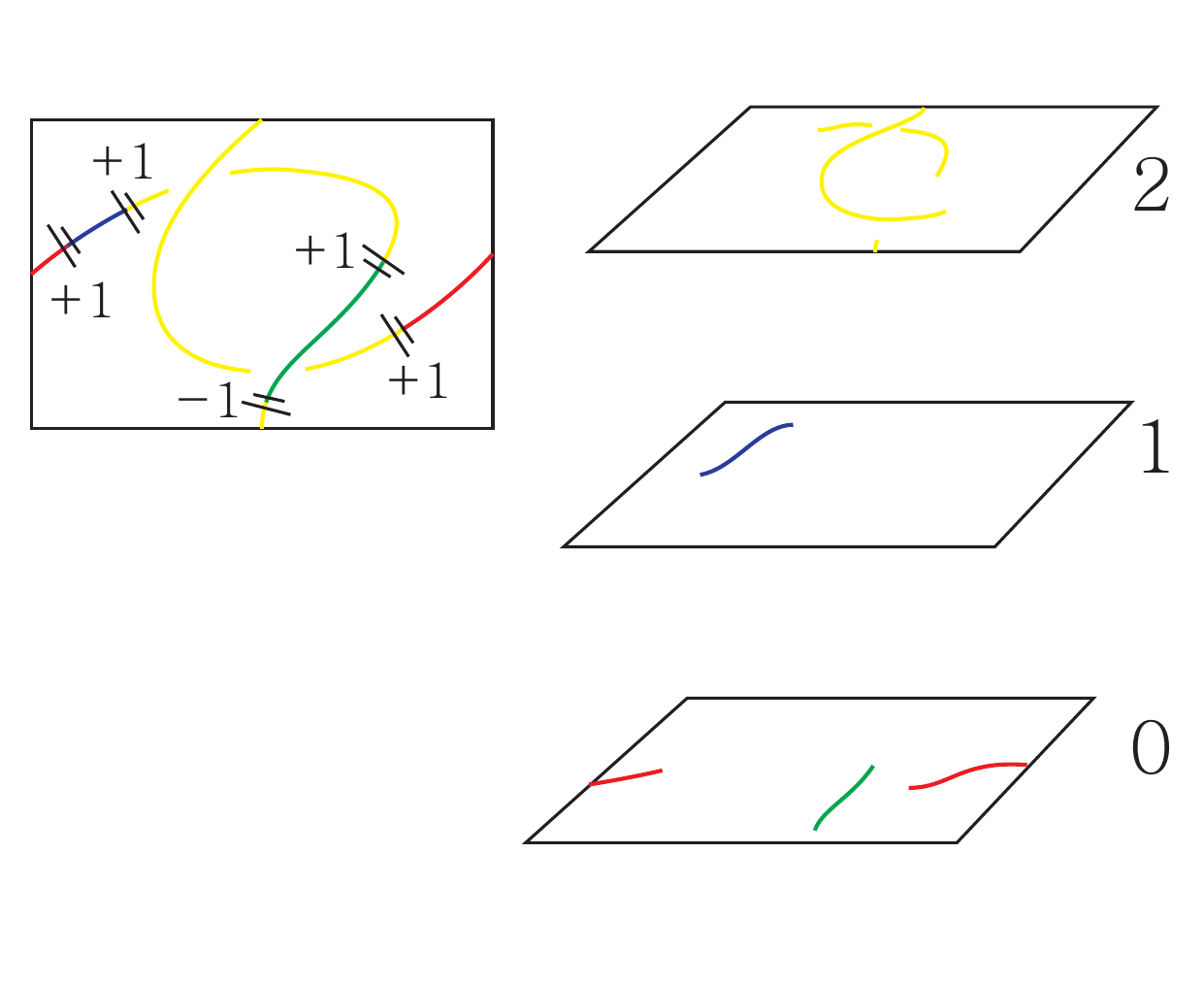}
\end{center}
\caption{Step 2 for knots of degree $3$}\label{fig:algo_lifting_deg3-2new}
\end{figure}

\textbf{Step 2.} Let us imagine $k$ parallel planes placed as Fig.~\ref{fig:algo_lifting_deg3-2new}. Give numbers to planes from bottom to top by integers from $0$ to $k-1$. Draw $k$ copies of $D$ for each plane. For a copy of a diagram on the plane with a number $s \in 
\mathbb{Z}_{k}$ erase long arcs which have not the height $s$. \\

\begin{figure}[h!]
\begin{center}
 \includegraphics[width = 8cm]{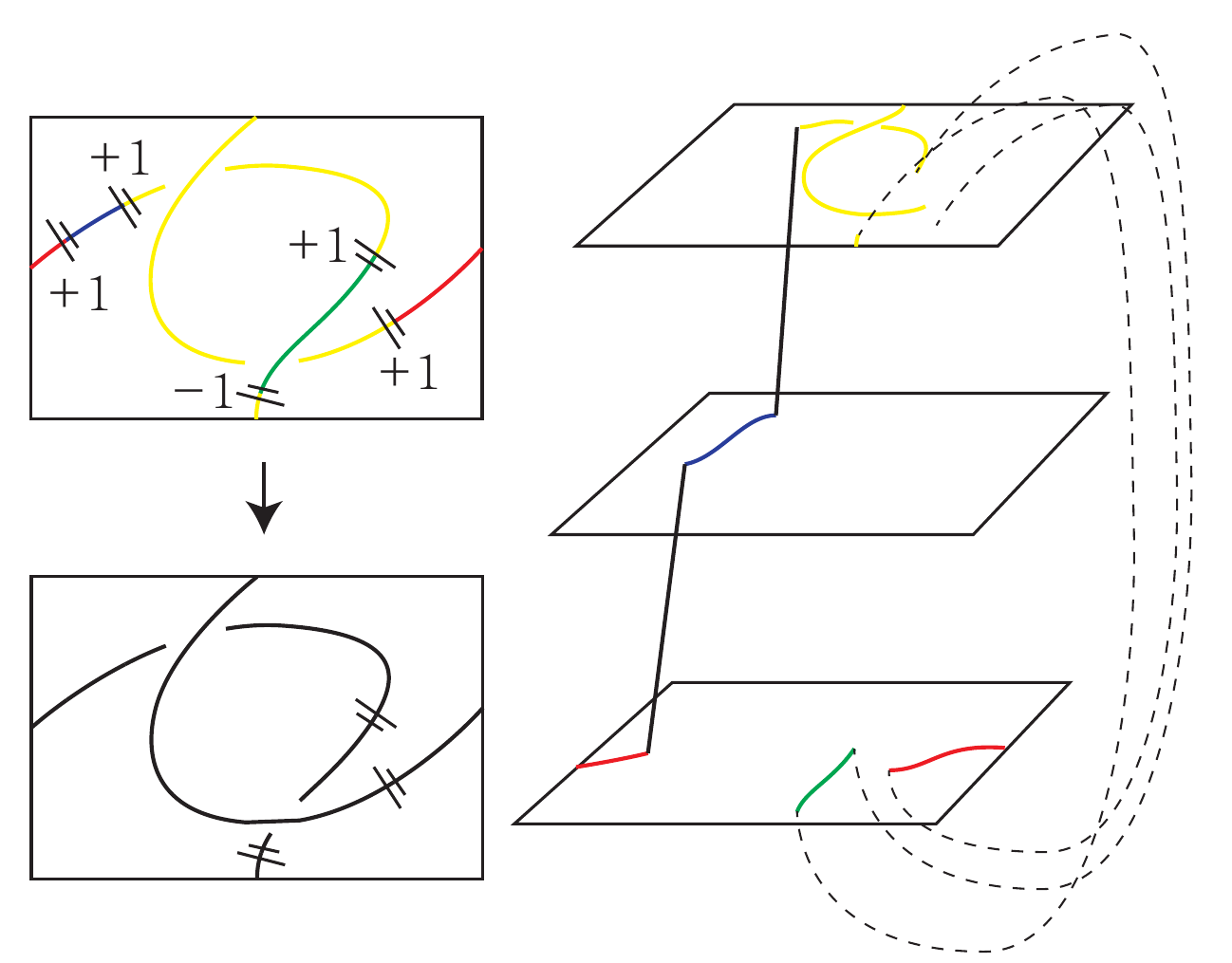}
\end{center}
\caption{Step 3 for knots of degree $3$}\label{fig:algo_lifting_deg3-3}
\end{figure}

\textbf{Step 3.} 
Let us suppose that we are walking on a long arc of $D$ with height $s \neq 0, k-1$. When we meet the double line, if it is a longer line, then we connect the corresponding arc on the plane with number $s$ to the arc on the plane with number $s+1$. If it is a shorter line, then we connect the corresponding arc on the plane with number $s$ to the arc on the plane with number $s-1$.
If the height is $s=k-1$ and we meet the longer line, then we connect the corresponding arc on the plane with number $k-1$ to the arc on the plane with number $0$ outside plane as described in Fig.~\ref{fig:algo_lifting_deg3-3}.
Similarly, if $s=0$ and we meet the shorter line, then we connect to the act on the plane with number $k-1$ outside plane.

From the above visualization of $\hat{K}$ we can obtain a diagram $\hat{D}$ from $\hat{K}$ as follows: for a diagram $D$ with double lines, we give heights to long arcs of $D$. For each crossing, if two heights of over/under crossings are same, then we remain it. If two heights of over/under crossings are different, then we change (if it is necessary) over/under information so that the arc with bigger height becomes over crossing. We remove double lines except double lines connecting long arcs with heights $0$ and $k-1$. Let us denote the obtained diagram by $\hat{D}$, see Fig.~\ref{fig:algo_lifting_deg3-3}.


\begin{lem}\label{lem:lifting-degk}
Let $D$ and $D'$ be two oriented diagrams with double lines. If they are equivalent, then $\hat{D}$ and $\hat{D'}$ are equivalent as knots in $S_{g} \times S^{1}$ with degree $1$ or $-1$.
\end{lem}

\begin{proof}
    Assume that $D'$ is obtained from $D$ by applying one of the moves in Fig.~\ref{moves2}. If $D'$ is obtained from $D$ by applying moves (1), (2), (3), (1'), (2'), (3') and (3''), then it is easy to see that $\hat{D'}$ can be obtained from $\hat{D}$ by applying virtual Reidemeister moves.

Suppose that $D'$ is obtained from $D$ by applying the move (4). Assume that the long arcs in move (4) have heights $a,a+1$ and $b$. If $a$ and $a+1$ are different with $0$ and $k-1$, then the proof is same to the proof of Lemma~\ref{lem:lifting-deg0}. If $a=k-1$, then the double line in the move (4) is still placed in $\hat{D}$ and $\hat{D'}$ and they are equivalent by the move (4). Similarly, we can show that if $D'$ is obtained from $D$ by applying the move (5), then $\hat{D}$ and $\hat{D'}$ are equivalent and hence it completes the proof.
\end{proof}

\begin{figure}[h!]
\begin{center}
 \includegraphics[width = 8cm]{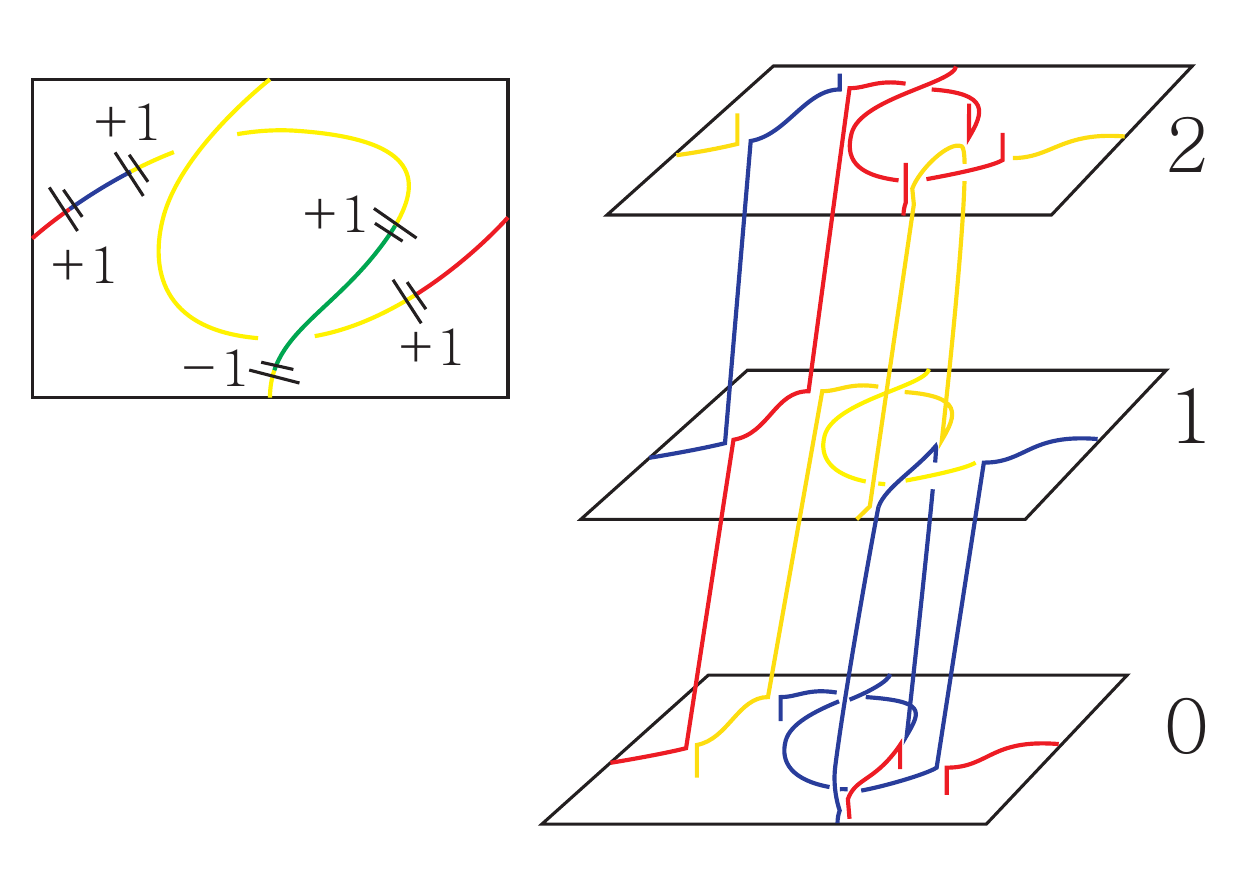}
\end{center}
\caption{Liftings $\{\hat{K}_{s}\}_{s \in \mathbb{Z}}$ of a knot $K$ of degree $k$ where $\hat{K}_{s}(0) = (x,e^{\frac{2\pi i}{k}{s}})$}\label{fig:exa_inf-lifting_deg3}
\end{figure}

If we consider liftings $\hat{K}_{s}$ such that $\hat{K}_{s}(0) = \hat{K}_{s}(0) = (x, e^\frac{2\pi si}{k})$, then we obtain a link in $S_{g} \times S^{1}$ of $k$ components, see Fig.~\ref{fig:exa_inf-lifting_deg3}. 

\begin{cor}\label{cor:multi-lifting-degk}
Let $K$ and $K'$ be knots in $S_{g} \times S^{1}$. If $K$ and $K'$ are equivalent, then the liftings $\sqcup_{s = 0}^{k-1} \hat{K}_{s}$ and $\sqcup_{s = 0}^{k-1} \hat{K}'_{s}$ are equivalent in $S_{g} \times S^{1}$.  
\end{cor}

\begin{figure}[h!]
\begin{center}
 \includegraphics[width = 8cm]{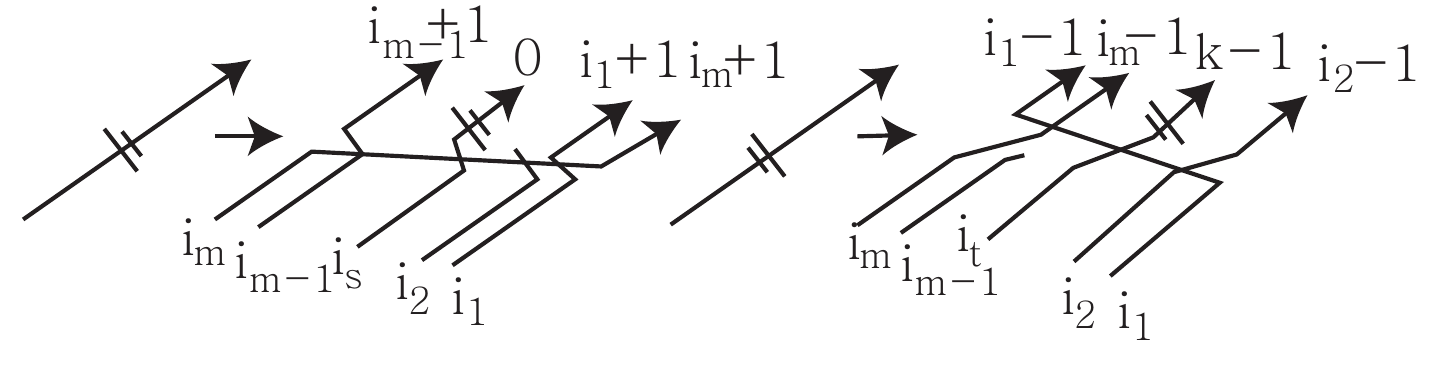}

\end{center}
 \caption{We number arcs as described in the present figure, but we place a double line, when $i_{s}=k-1$ and $i_{t}=0$, respectively.}\label{diag-const-lifting-numbering-degk}
\end{figure}

{\bf The construction of diagram of $\sqcup_{s = 0}^{k-1} K_{s}$ for a knot $K$ in $S_{g} \times S^{1}$ with degree $k$} is same to the case of knots $K$ in $S_{g} \times S^{1}$ with degree $0$ except Step 4. We replace Step 4 by Step 4':\\
{\bf Step 4'.} 
Walking along diagrams, we number arcs so that the number of an arc is not changed when we pass through the crossings corresponding to a classical crossing of a knot $D$. When we pass classical crossings associated to double lines, the numbering of arcs before twisting is changed by $\pm 1$ modulo $k$ with the rule described in Fig.~\ref{diag-const-lifting-numbering-degk}, but the numbering of arcs with crossings is same to the numbering before twisting. We place a double line on the arc with number $0$ or $k-1$ as described in Fig.~\ref{diag-const-lifting-numbering-degk}.

\begin{figure}[h!]
\begin{center}
 \includegraphics[width = 12cm]{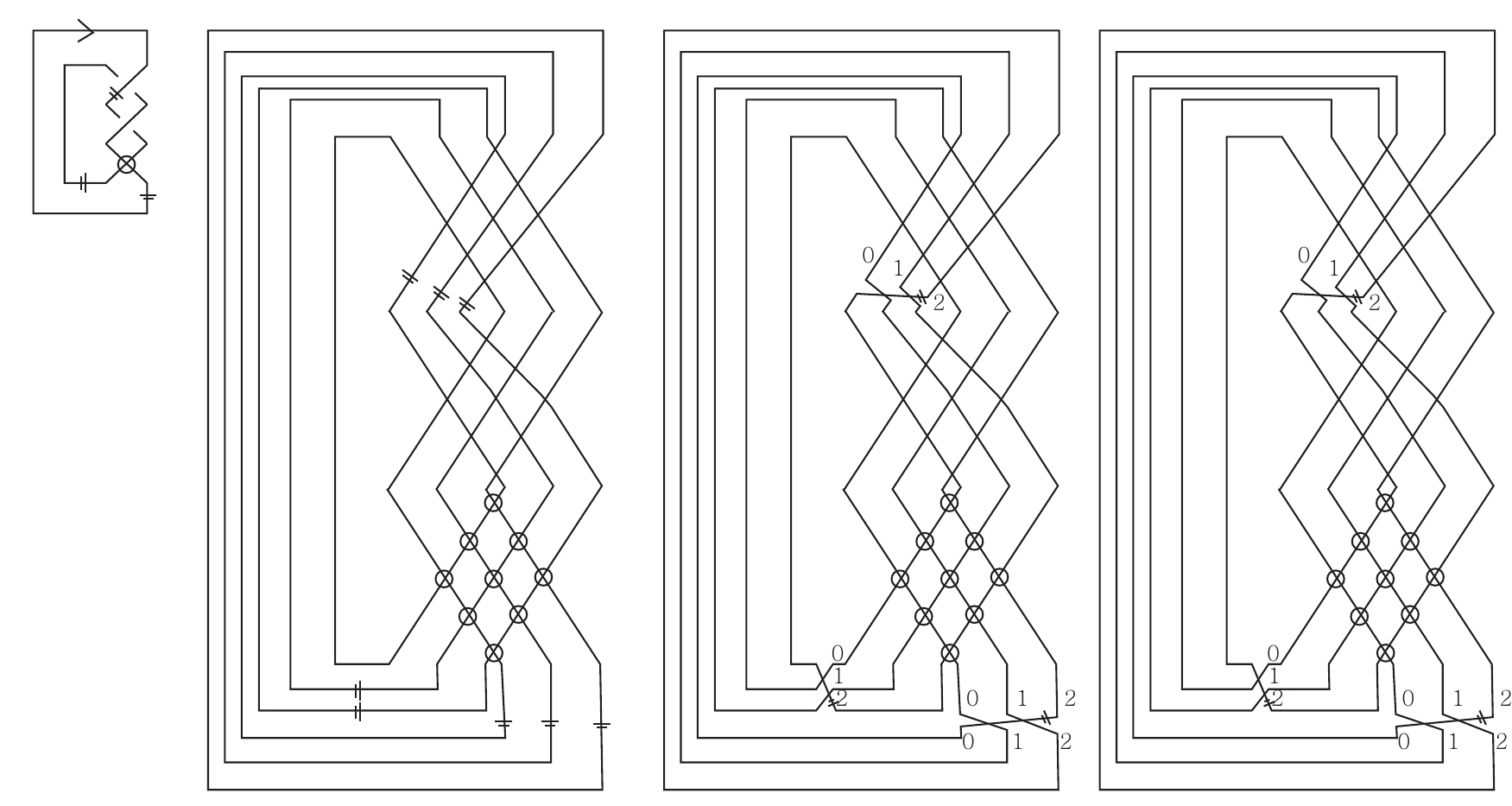}

\end{center}
 \caption{Construction of the diagram of 3-fold lifting of the knot in $S_{g} \times S^{1}$ of degree $3$}\label{exa-construction-lift-dl-degree3}
\end{figure}

We denote the obtained diagram by $\sqcup_{s=0}^{k-1} \hat{D}_{s}$, for example, see Fig.~\ref{exa-construction-lift-dl-degree3}.

\begin{proof}[Proof of Corollary \ref{cor:multi-lifting-degk}]
Corollary \ref{cor:multi-lifting-degk} also can be proved by using diagrams constructed as above. Let $D$ and $D'$ be diagrams with double lines associated to knots $K$ and $K'$ in $S_{g} \times S^{1}$ of degree $k$. If $D'$ is obtained from $D$ by applying moves (1), (2), (3), (1'), (2'), (3') and (3''), then it is easy to see that $\sqcup_{s=0}^{k-1} \hat{D'}_{s}$ can be obtained from $\sqcup_{s=0}^{k-1}\hat{D}_{s}$ as links in $S_{g}\times S^{1}$.

Suppose that $D'$ is obtained from $D$ by the move (4). Then $\sqcup_{s=0}^{k-1}\hat{D}_{s}$ and $\sqcup_{s=0}^{k-1}\hat{D'}_{s}$ have diagrams as described in Fig.~\ref{fig:multi-lifting-degk-move4}. 
\begin{figure}[h!]
\begin{center}
\includegraphics[width = 8cm]{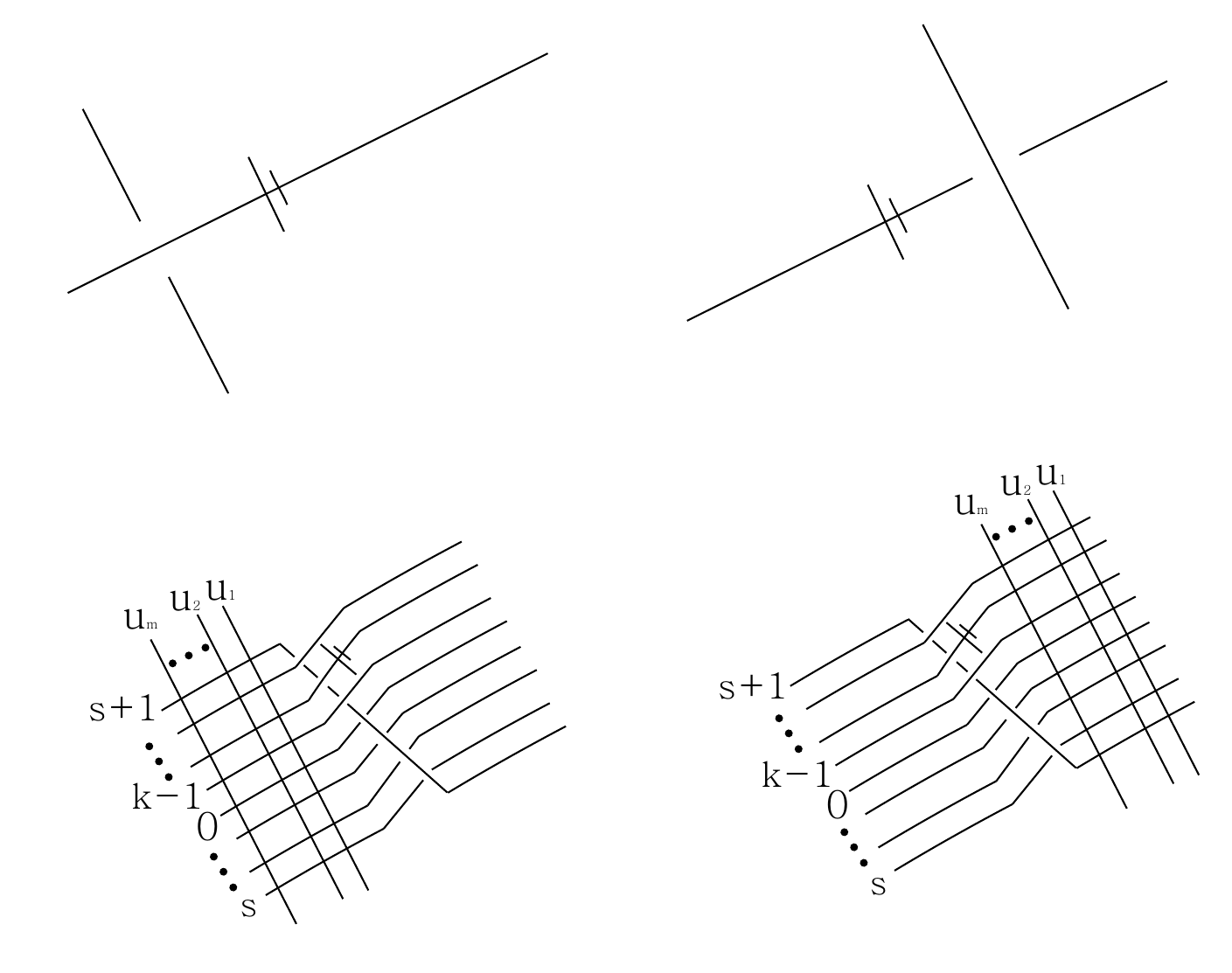}
\end{center}
\caption{Diagrams $\sqcup_{s=0}^{k-1}\hat{D}_{s}$ and $\sqcup_{s=0}^{k-1}\hat{D'}_{s}$ obtained from $D$ and $D'$ when $D'$ is obtained from $D$ by the move (4)}\label{fig:multi-lifting-degk-move4}
\end{figure}
It is enough to show that one arc with number $u_{i}$ can pass crossings. If $0\leq u \leq s$, then the arc with number $u$ can pass through the crossings by Reidemeister moves and the move (4) as described in Fig. \ref{fig:multi-lifting-degk-move4-case1}.
\begin{figure}[h!]
\begin{center}
\includegraphics[width = 10cm]{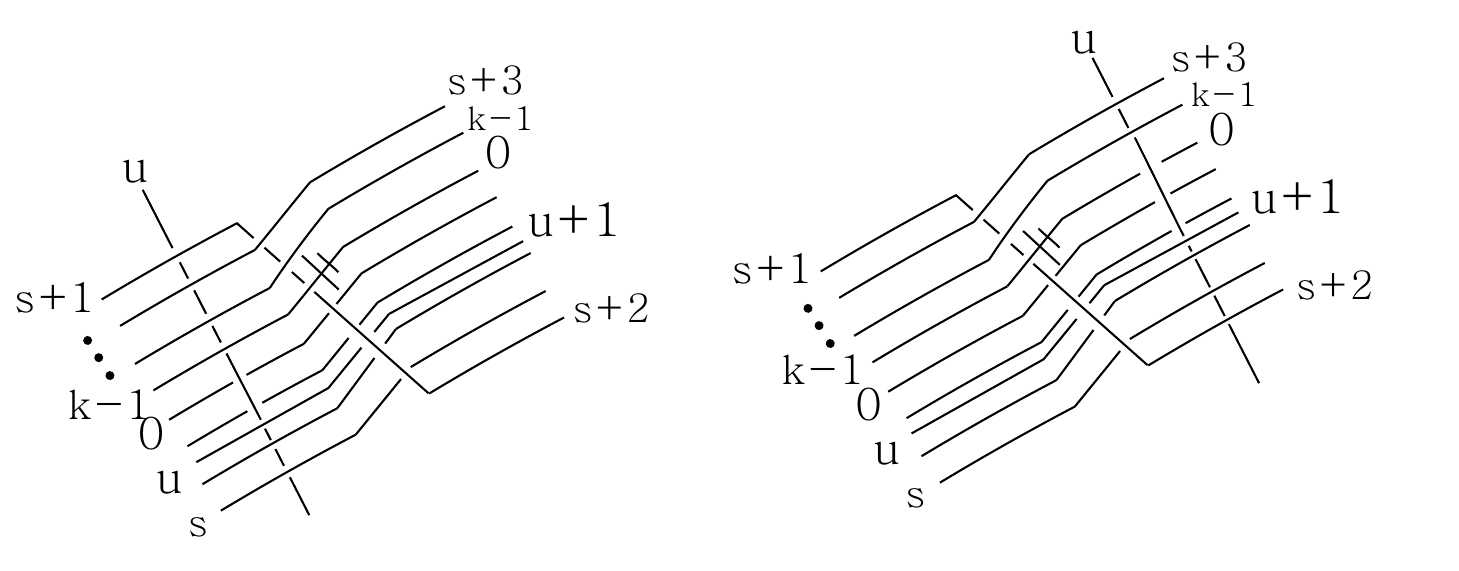}
\end{center}
\caption{Case 1. $0\leq u \leq s$}\label{fig:multi-lifting-degk-move4-case1}
\end{figure}
If $s+1\leq u \leq k-1$, then the arc with number $u$ can pass through the crossings by Reidemeister moves and the move (4) as described in Fig. \ref{fig:multi-lifting-degk-move4-case2}.
\begin{figure}[h!]
\begin{center}
\includegraphics[width = 10cm]{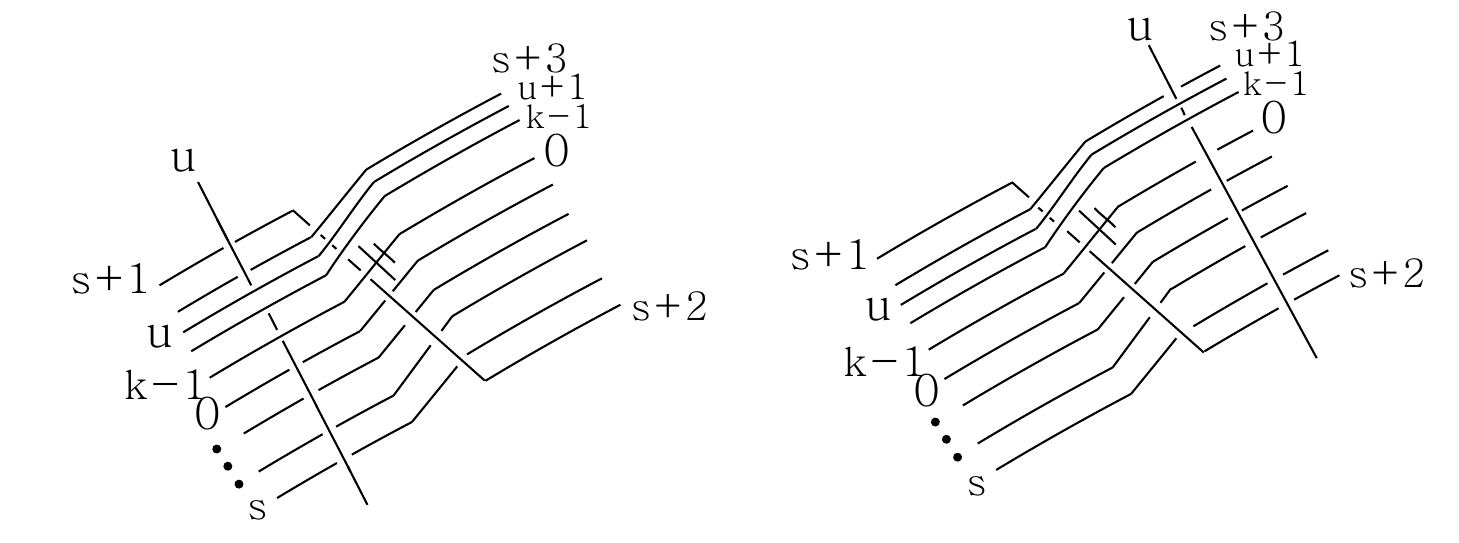}
\end{center}
\caption{Case 2. $s+1\leq u \leq k-1$}\label{fig:multi-lifting-degk-move4-case2}
\end{figure}

If $D'$ is obtained from $D$ by the move (5), then $\sqcup_{s=0}^{k-1}\hat{D'}_{s}$ can be obtained from $\sqcup_{s=0}^{k-1}\hat{D}_{s}$ by a series of the second Reidemeister moves and the move (5). Let us leave the details for the case of the move (5) to readers.

\end{proof}

\section{Liftings of knots in $S_{g} \times S^{1}$ and cut system}

\begin{figure}[h!]
\begin{center}
 \includegraphics[width = 4cm]{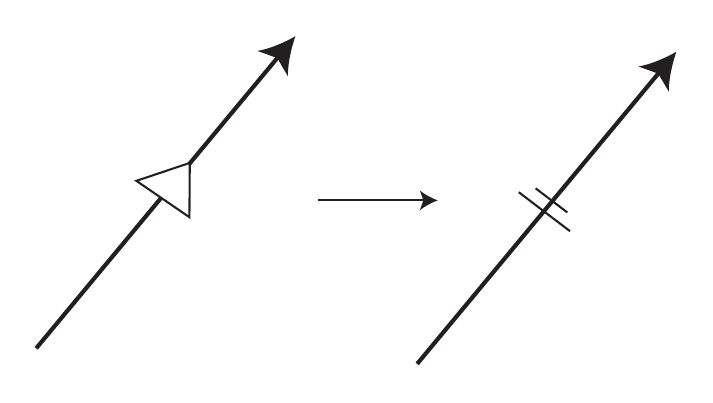}
\end{center}
\caption{Map from $(D,P)$ to a knot in $S_{g} \times S^{1}$}\label{fig:map-cutpt-doubleline}
\end{figure}

For a virtual link diagram with a cut system $(D,P)$, one can obtain a link diagram with double lines $D^{dl}$ by replacing cut points to double lines as described in Fig.~\ref{fig:map-cutpt-doubleline}. 

\begin{lem}
For two virtual link diagrams $(D_{1},P_{1})$ and $(D_{2},P_{2})$ with cut systems, if $(D_{1},P_{1})$ and $(D_{2},P_{2})$ are equivalent under oriented cut point moves, then the link diagrams $D_{1}^{dl}$ and $D_{2}^{dl}$ with double lines obtained from $(D_{1},P_{1})$ and $(D_{2},P_{2})$ respectively, are equivalent.
\end{lem}

\begin{proof}
Suppose that $(D_{2},P_{2})$ is obtained from $(D_{1},P_{1})$ by one of oriented cut point moves. It is easy to see that if $(D_{2},P_{2})$ is obtained from $(D_{1},P_{1})$ by the move of ``passing'' of a cut point through virtual crossing or the cancellation move of an adjacent coherent and incoherent cut points, then $D_{1}^{dl}$ and $D_{2}^{dl}$ are equivalent as knots in $S_{g} \times S^{1}$. If $(D_{2},P_{2})$ is obtained from $(D_{1},P_{1})$ by deletion of four cut points around a classical crossing, then one can show that $D_{1}^{dl}$ and $D_{2}^{dl}$ are equivalent as knots in $S_{g} \times S^{1}$ as described in Fig.~\ref{fig:proof-well-defined-map}.
 \begin{figure}[h!]
\begin{center}
 \includegraphics[width = 8cm]{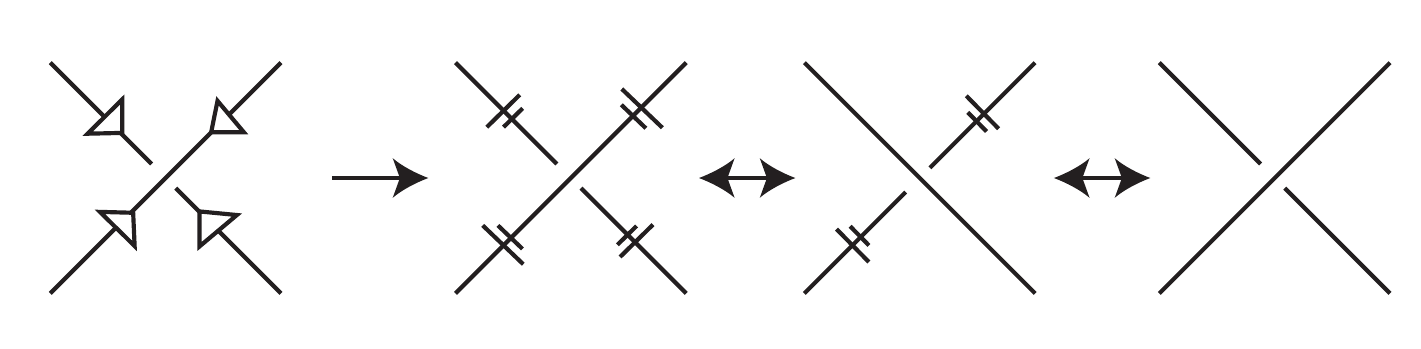}
\end{center}
\caption{Equivalence relation in virtual link diagrams with cut systems is preserved}\label{fig:proof-well-defined-map}
\end{figure}
\end{proof}
\begin{cor}
For a virtual link diagram $(D,P)$ with a cut system, the obtained link diagram $D^{dl}$ with double lines has degree $0$. 
\end{cor}

\begin{proof}
By definition of cut system it is obvious that the numbers of coherent and incoherent cut points are same. Therefore, the obtained link diagram $D^{dl}$ with double lines has degree $0$.
\end{proof}

\begin{lem}
Let $D$ be an oriented diagram of a knot in $S_{g} \times S^{1}$ of degree $0$. Then the $m$-fold lifting of $D$ is mod $m$ almost classical for any $m >1$.
\end{lem}
\begin{proof}
Let $\hat{D}^{m}:= \sqcup_{s=1}^{m}\hat{D}_{s}$ be an oriented diagram of a $m$-fold lifting of $D$. Then $\hat{D}^{m}$ has $m$ parallel of $D$ with twists in the places of double lines. Let us color it by $\mathbb{Z}_{m}$. Fix any point on $D$ and color $m$ semi-arcs of $\hat{D}^{m}$ corresponding to the fixed point by $0,1, \dots, m-1$ from right to left. For convenience let us denote the coloring by $(0,1,\dots, m-1)$.

\noindent {\bf Step 1.} In the diagram $\tilde{D}^{m}$ there are two kinds of twists of arcs and the arcs near to the twists are colored as described in Fig.~\ref{fig:proof-step1}. That is, the order of colors of $m$ parallels remains. Therefore, it suffices to show that $m$ parallel of $D$ is mod $m$ Alexander colorable.

\begin{figure}[h!]
\begin{center}
 \includegraphics[width = 6cm]{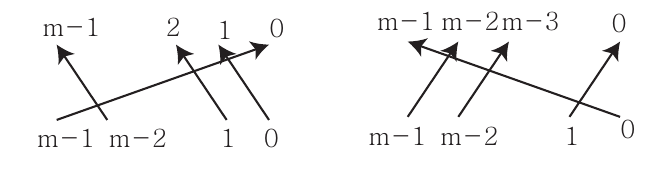}
\end{center}
\caption{Proof step 1}\label{fig:proof-step1}
\end{figure}

\noindent {\bf Step 2.} Each crossing of $D$ corresponds to $m^{2}$ crossings as described in Fig.~\ref{fig:proof-step2-1}. Note that $i$-th arcs from right of the parallel corresponding to two arcs for a crossing of $D$, are in same component. 
\begin{figure}[h!]
\begin{center}
 \includegraphics[width = 6cm]{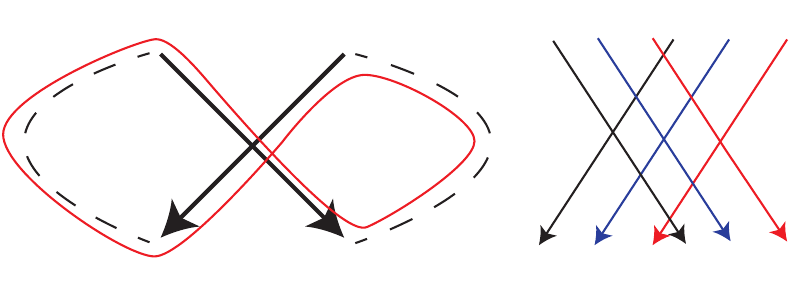}
\end{center}
\caption{Proof order of parallel}\label{fig:proof-step2-1}
\end{figure}
Now $m$ parallels going into the part of $m^{2}$ crossings are colored by $(0,1,\dots, m-1)$. Then in $\mathbb{Z}_{m}$ the Alexander numbering is preserved passing through $m^{2}$ crossings.
\begin{figure}[h!]
\begin{center}
 \includegraphics[width = 6cm]{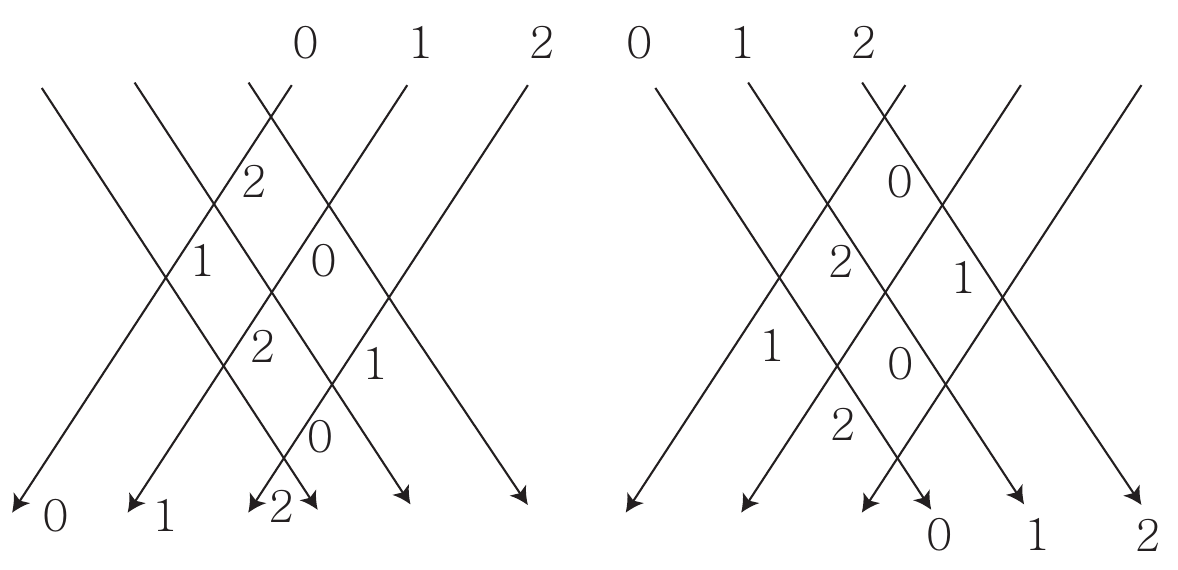}
\end{center}
\caption{Proof step 3}\label{fig:proof-step2-2}
\end{figure}

\noindent {\bf Step 3.} One can show that the numberings for arcs corresponding to $m^{2}$ crossings with two entering colorings $(0,1,\dots, m-1)$ satisfy the properties in Fig.~\ref{fig:proof-step3}.

\begin{figure}[h!]
\begin{center}
 \includegraphics[width = 3cm]{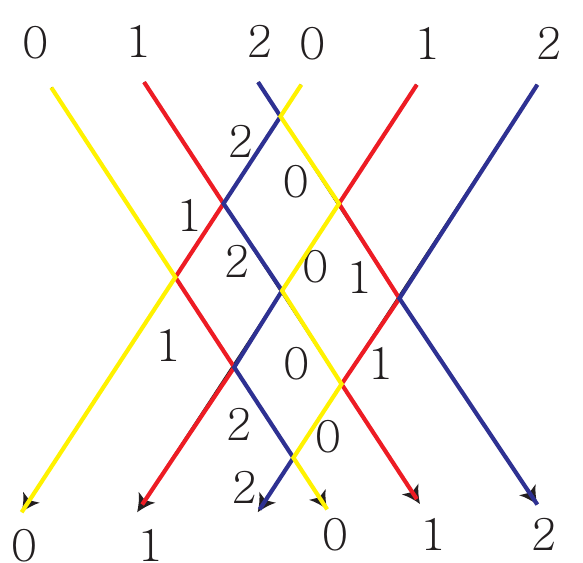}
\end{center}
\caption{Proof step 4}\label{fig:proof-step3}
\end{figure}
From the previous steps, the proof is completed.
\end{proof}

\begin{thm}
Let $D^{dl}$ be a link diagram with double lines obtained from a virtual link diagram $(D,P)$ with cut system. Then the $m$-fold lifting of $D^{dl}$ is mod $m$ almost classical for any $m >1$.
\end{thm}

\begin{exa}
We obtain $3$-fold covering of $D^{dl}(3_{1})$ from the trefoil knot $3_{1}$ and it must be mod $3$ almost classical. In Fig.~\ref{fig:exa-Alexander-num-SgS1} one can find a mod $3$ Alexander numbering.

\begin{figure}[h!]
\begin{center}
 \includegraphics[width = 6cm]{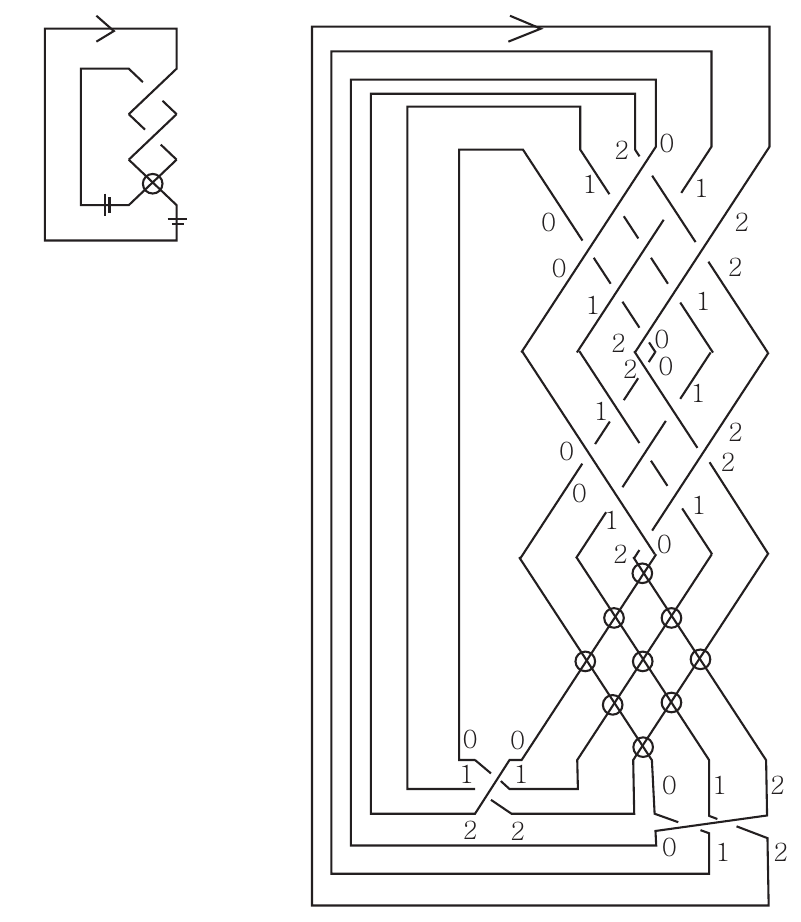}
\end{center}
\caption{Example of Alexander numbering of $D^{dl}$}\label{fig:exa-Alexander-num-SgS1}
\end{figure}

\end{exa}

\begin{rem}
We have two liftings of virtual knots, which provide mod $m$ almost classical links. The difference is the following: In \cite{NaokoKamada} the constructed lifting can be obtained from $m$ parallels of virtual knot such that classical crossings between different components are replaced by virtual crossings. But in our construction, classical crossings between different crossing remains as classical crossings.
\end{rem}

\end{document}